\documentclass[preprint,12pt,authoryear]{elsarticle}
\makeatletter
\IfFileExists{TRC_revision/main.tex}{\def\rootdir{TRC_revision/}}{\def\rootdir{}}
\edef\input@path{{\rootdir}}
\makeatother

\usepackage[margin=1in]{geometry}
\usepackage{amsmath,amssymb,amsthm,mathtools}
\usepackage{booktabs,tabularx,multirow,array}
\usepackage{threeparttable}
\usepackage{graphicx,subcaption,placeins}
\edef\setgpath{\noexpand\graphicspath{{\rootdir}{\rootdir fig/}}}\setgpath
\usepackage[ruled,vlined,linesnumbered]{algorithm2e}
\usepackage{enumitem}
\usepackage{xcolor}
\usepackage{xurl}
\usepackage[colorlinks=true,allcolors=blue]{hyperref}
\usepackage[nameinlink]{cleveref}
\usepackage{tikz}
\usetikzlibrary{arrows.meta,positioning,calc}
\usepackage{siunitx}

\newtheorem{proposition}{Proposition}

\newcommand{\sV}{\mathcal{V}}
\newcommand{\sC}{\mathcal{C}}
\newcommand{\sL}{\mathcal{L}}

\newcommand{\sN}{\mathcal{N}}
\newcommand{\bigM}{\mathbb{M}}
\crefname{equation}{}{}

\Crefname{equation}{Eq.}{Eqs.}
\crefname{algocf}{Algorithm}{Algorithms}
\crefname{algorithm}{Algorithm}{Algorithms}
\Crefname{algorithm}{Algorithm}{Algorithms}
\crefname{section}{Section}{Sections}
\crefname{subsection}{Section}{Sections}
\crefname{figure}{Fig.}{Figs.}
\crefname{table}{Table}{Tables}
\crefname{proposition}{Proposition}{Propositions}

\journal{Transportation Research Part C: Emerging Technologies}

\begin{document}

\begin{frontmatter}

\title{National-scale scheduling of long-distance electric vehicle charging under plug capacity}

\author[tamu,anl]{Taner Cokyasar}
\ead{tcokyasar@tamu.edu}
\affiliation[tamu]{organization={Department of Engineering Technology and Industrial Distribution, Texas A\&M University},
            addressline={3367 TAMU},
            city={College Station},
            postcode={77843},
            state={TX},
            country={USA}}
\affiliation[anl]{organization={Argonne National Laboratory},
            addressline={9700 S. Cass Ave.},
            city={Lemont},
            postcode={60439},
            state={IL},
            country={USA}}

\begin{abstract}
Can today's U.S. public charging network serve long-distance trips by electric vehicles (EVs), and at what cost to travelers? We plan the charging stops of 2.2 million long-distance trips that need en-route charging, from a 1\% sample of a one-month simulation, through the contiguous U.S. fast-charging network. Plans value detour, charging, and waiting time. Charging to a fixed battery level makes each trip's problem a shortest path, solved exactly; common modeling shortcuts can reverse the cost comparison. For trips competing for plugs, we compare four coordination methods, from none to reservations by departure time; the best two come within 2.2\% and 6.1\% of the best known solutions on small cases. First, if every plug serves every EV, existing fast chargers can reach 98.5\% of the trips that need to charge, excluding those with an overly long link, with under four stops and about seven minutes of detour. Second, charging speed and electricity price decide the cost comparison: at \$0.15 per kWh, the EV is cheaper only for travelers valuing their charging time below \$11 per hour with 50~kW, or \$111 per hour with 350~kW chargers; at \$0.45 per kWh, it costs more at either speed. Third, without coordination, long queues form at several hundred small stations at 50~kW while nearby larger ones are less busy. If every trip books a plug, almost all waiting disappears on all trips at 350~kW for under a minute of extra detour, and at \$0.15 per kWh the EV trip stays about 21\% cheaper.
\end{abstract}


\end{frontmatter}

\section{Introduction}\label{sec:intro}

Long-distance travel is where battery electric vehicles (EVs) differ most from internal combustion engine (ICE) vehicles. A trip of several hundred miles requires one or more en-route charging stops; each stop needs a detour to a station, a charging session that lasts tens of minutes at direct-current fast chargers (DCFC) and hours at Level~2 (L2) chargers, and possibly a wait for a free plug. Whether the existing public network can absorb this demand, and at what cost to travelers, depends jointly on where stations are, how many plugs they have, how fast they charge, and how many EVs try to use them at the same time.

National assessments of this question have so far relied on coverage metrics or aggregate demand models. \citet{hanig_2025} measure how much of the U.S.\ road network is consecutively reachable from existing fast chargers; \citet{wood_2023} estimate how many plugs a 2030 fleet will need from simulated charging demand; \citet{borlaug_2025} evaluate the economics of corridor fast-charging stations. These studies answer where coverage is missing or how much infrastructure to build, but they do not schedule individual trips through the existing network, so they cannot say how individual travelers are affected when demand concentrates on particular stations at particular times. Optimization models that do schedule individual vehicles~\citep{schneider_2014,keskin_2019,kullman_2021,dastpak_2023} are solved for one vehicle or a small fleet; models that capture station congestion at the network level~\citep{wang_2019,chen_2020,kinay_2023} are siting models with stylized demand.

This paper schedules the charging stops of the electrified long-distance trips of a national sample that need en-route charging through the real station inventory. It takes as given (i) about 2.8 million long-distance personal-vehicle trips from a 1\% national sample produced by the POLARIS agent-based model~\citep{auld_2016}, each routed link by link on the national road network, and (ii) the 14,260 public DCFC stations (68,362 plugs) and 64,477 public L2 stations (169,385 plugs) in the Alternative Fuels Station Locator~\citep{afdc_2026}. For each trip it chooses after which links to stop and at which nearby station, so as to minimize the value of detour time, charging time, and waiting time, subject to state-of-charge (SoC) limits and plug capacity. It then evaluates the outcome with a three-tier comparison against an ICE vehicle making the same trip.

The paper makes four contributions.
\begin{enumerate}[leftmargin=*]
\item \emph{A per-vehicle formulation that is exact under charge-to-target.} Under the common policy of charging to a fixed target SoC, the capacity-free problem of one vehicle is a shortest path on a \emph{charging-opportunity} directed acyclic graph (DAG) whose nodes are (link endpoint, candidate station) pairs. Because the SoC after every stop is known in advance, SoC needs no discretization and the dynamic program (DP) is exact. It schedules 2.2 million national trajectories with three candidate stations per link endpoint without a commercial solver. The DP is a standard shortest path once the policy is fixed; its value here is that it makes the national evaluation possible.
\item \emph{An evaluation of congestion-resolution heuristics against an exact reservation-type model.} We formulate the multi-vehicle problem as a continuous-time mixed-integer program (MIP) with per-plug sequencing and compare it, on small corridor instances, with four ways of coupling the per-vehicle DPs: no coordination (the DPs followed by a first-in, first-out (FIFO) queue), an iterative station-penalty heuristic, a FIFO best-response heuristic, and sequential reservation, in which trips plan in order of departure against the plug-times booked by earlier trips. Its time-dependent DP against per-plug calendars, with merged busy blocks and a fast path for trips that need not wait, makes a national reservation pass feasible at 350~kW. All heuristic schedules are scored with the same FIFO simulation, in which a vehicle's queueing delay propagates to its later stops.
\item \emph{A quantification of modeling choices that matter at scale.} A stepwise comparison on the full national data set shows how much each modeling choice changes the answer: how the starting SoC is bounded, whether SoC is discretized, and whether charging time is valued. Each of these choices changes feasibility, the number of stops, or the cost comparison by a large factor.
\item \emph{A national stress test of the existing network, uncoordinated and coordinated.} We report feasibility, detours, charging time, queueing, and EV-versus-ICE costs for DCFC at 50 and 350~kW and L2 at 5 and 15~kW; how early uncoordinated waiting appears as demand grows; how much of it sequential reservation removes (all trips at 350~kW, a 10\% sample at 50~kW); and sensitivity analyses on EV range, energy use, SoC limits, starting SoC, charger power, candidate-set size, detour circuity, detour speed, and plug supply.
\end{enumerate}

The rest of the paper proceeds as follows. \Cref{sec:lit} reviews studies on where to build charging stations, how individual vehicles plan their charging stops, how many vehicles can be coordinated at shared stations, and how national charging needs have been assessed. \Cref{sec:problem} states the problem, namely which stations each trip uses and when, given battery limits and the number of plugs at each station, and writes it as an optimization model. \Cref{sec:method} explains how we solve it at national scale: an exact method that plans one vehicle at a time, and four ways to handle the competition among vehicles for the same plugs. \Cref{sec:design} describes the national data and the experiments. \Cref{sec:validation} tests the methods against the exact model on small corridors and shows how much common modeling shortcuts change the results. \Cref{sec:national} applies the framework to the national data: whether the existing network can serve the trips, what the trips cost compared with gasoline vehicles, where and when queues form, and how much plug reservations reduce them. \Cref{sec:discussion} discusses what the results imply and what they do not, states the limitations, and concludes.
\section{Related work}\label{sec:lit}

Four lines of work bear on this paper: the design of charging networks for long-distance travel, en-route charging decisions for individual vehicles, coordination of many vehicles at capacitated stations, and national assessments of charging infrastructure and travel cost. \Cref{tab:lit} positions representative studies.

Charging networks for long-distance travel are usually designed with flow-based location models, which site stations so that trips can be completed within a driving range. The flow-refueling location model~\citep{kuby_2005} and its capacitated extension~\citep{upchurch_2009} treat station locations and sizes as decisions and origin--destination flows as data; multi-period deployment~\citep{li_2016}, U.S.-wide placement for long-distance trip completion~\citep{he_2019}, and joint siting for local and long-distance travel with adoption feedback~\citep{anjos_2020} follow the same logic. Several models represent queueing at stations explicitly: expanded charging networks with an approximate queueing-time function inside a variational-inequality equilibrium~\citep{wang_2019}, bilevel location and capacity design with route-choice and waiting-time equilibrium~\citep{chen_2020}, siting in congested intercity networks with detours~\citep{bao_2021,zeng_2024}, profit-maximizing station design with elastic demand and congestion~\citep{huang_2020}, and bilevel location and sizing under a probabilistic waiting-time requirement~\citep{kinay_2023}. These models optimize the supply side and represent demand as aggregate flows on stylized or regional networks. We do the opposite: supply is fixed at the real station inventory, and demand is a large set of individual trajectories, so the question is how well the given network serves the given trips.

For a single vehicle on a known route, choosing where to charge is a resource-constrained shortest-path problem with SoC as the resource~\citep{irnich_2005}; \citet{baum_2015} compute such shortest feasible paths with charging stops on road networks. Electric vehicle routing problems add charging-station visits to vehicle routing with time windows~\citep{schneider_2014}, time-dependent waiting times at stations~\citep{keskin_2019}, and nonlinear charging functions with partial charging~\citep{montoya_2017}; optimal and adaptive recharging policies for a single EV on a network with uncertain station availability are characterized by \citet{sweda_2017}. For travelers using public stations, routing policies that anticipate stochastic station queues~\citep{kullman_2021}, adaptive multi-criteria route planning that trades driving, waiting, and charging time~\citep{schoenberg_2021}, and dynamic policies that react to occupancy information~\citep{dastpak_2023} have been proposed, and route planning through the real German charging network has been used to evaluate that network for typical EVs on 60 routes~\citep{hecht_2021}. These methods treat other vehicles as exogenous: queues are random or observed, not caused by the vehicles being planned. \citet{mandal_2025} solve the fixed-route charging problem of one vehicle exactly by dynamic programming with charging to a few discrete levels. Our per-vehicle DP uses the same shortest-path insight but exploits the charge-to-target policy to remove the SoC dimension, which is what makes it fast enough to solve millions of trajectories; the price is that partial charging and charging curves, which those models represent, are approximated here (\cref{sec:discussion}).

When many vehicles share plugs, individually optimal plans can overload stations. Highway studies coordinate charging schedules through peer-to-peer reservations~\citep{delrazo_2016} or distributed control of station queues~\citep{gusrialdi_2017}; reservation algorithms assign charging slots at stations~\citep{flocea_2022}, and one such scheme reduces queueing at urban chargers by about 70\% on the Sioux Falls network~\citep{lin_2025}; transit scheduling with capacitated chargers has been solved by column generation~\citep{devos_2024}. These approaches either handle one corridor or one fleet, or rely on a central operator that can assign charging slots. Our setting has millions of independent travelers on a national network. We therefore examine what can be achieved with simple coordination schemes that re-use the per-vehicle DP, and we measure how far they are from an exact model on instances small enough for it.

National studies have quantified the potential for electrifying personal travel from travel-survey data~\citep{needell_2016}, the gaps in consecutive fast-charging coverage of U.S.\ roads~\citep{hanig_2025}, the plugs needed for a 2030 fleet~\citep{wood_2023}, the levelized cost of charging~\citep{borlaug_2020}, the economics of corridor fast-charging stations~\citep{borlaug_2025}, observed fast-charger use~\citep{neaimeh_2017}, and the total cost of ownership of EVs~\citep{hagman_2016}. None schedules individual trips through the existing stations with plug capacity, so none reports the detours, charging times, and queues that a given electrification level produces on today's network, or how those translate into travel cost.

\begin{table}[!htb]
\centering
\caption{Positioning of representative studies.}\label{tab:lit}
\footnotesize
\setlength{\tabcolsep}{4pt}
\begin{tabularx}{\textwidth}{>{\raggedright\arraybackslash}p{4.3cm} ccc >{\raggedright\arraybackslash}X >{\raggedright\arraybackslash}X}
\toprule
Study & Siting & Scheduling & Capacity & Demand & Scale and method \\
\midrule
\citealp{kuby_2005,upchurch_2009} & \checkmark & & \checkmark$^a$ & O--D flows & Regional; MIP \\
\citealp{he_2019,anjos_2020} & \checkmark & & & O--D flows & National/regional; MIP \\
\citealp{wang_2019,chen_2020,kinay_2023} & \checkmark & & \checkmark & Equilibrium flows & Test networks; bilevel or VI \\
\citealp{bao_2021,zeng_2024} & \checkmark & & & O--D flows, detours & Intercity networks; exact/approx.\ partitioning \\
\citealp{schneider_2014,keskin_2019} & & \checkmark & \checkmark$^b$ & Fleet routes & Benchmark instances; MIP, metaheuristics \\
\citealp{kullman_2021,dastpak_2023} & & \checkmark & \checkmark$^c$ & One vehicle & Single trips; dynamic programming \\
\citealp{delrazo_2016,gusrialdi_2017} & & \checkmark & \checkmark & Highway corridor & One corridor; reservation, distributed control \\
\citealp{devos_2024} & & \checkmark & \checkmark & Transit blocks & One fleet; column generation \\
\citealp{hanig_2025,wood_2023} & & & & Trips/charging demand & National; coverage metric, simulation \\
\midrule
This paper & & \checkmark & \checkmark & 2.2M individual trajectories & National; exact per-vehicle DP and coordination heuristics benchmarked against an exact MIP on small corridors (11 of 25 proven optimal) \\
\bottomrule
\multicolumn{6}{p{0.97\textwidth}}{\emph{$^a$Capacitated variant~\citep{upchurch_2009}. $^b$Time-dependent waiting~\citep{keskin_2019}. $^c$Exogenous stochastic queues.}}\\
\end{tabularx}
\par\smallskip
\begin{minipage}{\linewidth}\footnotesize\textit{Note:} Siting: station locations or sizes are decisions. Scheduling: the timing and location of individual vehicles' charging stops are decisions. Capacity: plug or station capacity is modeled.\end{minipage}
\end{table}
\section{Problem definition}\label{sec:problem}

\subsection{Setting and assumptions}
Consider a set of electrified trips $\sV$. Trip $v \in \sV$ follows a fixed route $\sL_v = (1, \ldots, N_v)$ of directed road links obtained from a static shortest-path assignment, departs at a known time $T^0_v$ with a known starting SoC $B_v$, and travels at the link travel times of the assignment. A set of public charging stations $\sC$ is given; station $c \in \sC$ has $Z_c$ plugs, all of the scenario's charger power $W$. The decision for each trip is after which links to leave the route and at which station to charge; how long to charge follows from the target SoC of (A4) (\cref{fig:problem}). The objective is to minimize the total value of detour time, charging time, and waiting time over all trips, subject to SoC limits and plug capacity. The model rests on the following assumptions.

\begin{enumerate}[label=(A\arabic*),leftmargin=*,noitemsep,topsep=2pt]
\item Station locations and plug counts are fixed; this is a scheduling problem, not a siting or sizing problem.
\item All plugs in a scenario have the same power $W$; DCFC and L2 networks are studied as separate scenarios.
\item At the end of link $i$, a trip may detour to one of the $K$ stations nearest to the link endpoint, $\sC_{vi} \subseteq \sC$. In the national experiments the nearest station is always a candidate; the second to $K$-th nearest are candidates only if they are within $\bar{D}$ miles. Detour distance is the $L_1$ distance in projected coordinates multiplied by a circuity factor $\vartheta \geq 1$ and is driven at constant speed $\bar{s}$; the vehicle returns to the same link endpoint.
\item Every stop charges to the target SoC $\overline{S}$. The starting SoC may exceed $\overline{S}$ (most long-distance trips start after charging at home); $\overline{S}$ limits only en-route charging. SoC may never fall below $\underline{S}$.
\item Energy use is proportional to distance at a trip-specific rate $\rho_v$ (Wh/mile); charging power is constant, so the charging rate is $\phi_v = 10^5\, W/Q_v$ (\%/h) with battery capacity $Q_v$ (Wh). Charging to $\overline{S} = 80\%$ keeps sessions in the range in which constant power is a reasonable approximation.
\item Plugs serve arriving vehicles first-in, first-out; after departure, vehicles do not deviate from their plans: they do not balk, reroute, or change departure times when stations are busy. This rule governs the evaluation and the FIFO-based heuristics; the exact model of \cref{sec:mip} relaxes it and may sequence the stops on a plug in any order (\cref{sec:eval}). Each trip is independent; multi-day tours are not chained.
\end{enumerate}

\begin{figure}[!htb]
\centering
\begin{tikzpicture}[
    det/.style={->, >=stealth, dashed, thick, black},
    station/.style={draw, circle, minimum size=0.8cm, fill=black!6, font=\footnotesize, thick},
    sl/.style={font=\footnotesize}]
\def\lw{1.7}
\fill[black!7] (0,0) rectangle (7*\lw,0.5);
\draw[thick] (0,0) rectangle (7*\lw,0.5);
\foreach \k in {1,...,6}{\draw[black!40] (\k*\lw,0) -- (\k*\lw,0.5);}
\foreach \i in {1,...,7}{\node[font=\small] at ({(\i-0.5)*\lw},0.25) {$\ell_{\i}$};}
\draw[->,>=stealth,thick] (0.2,0.75) -- (7*\lw-0.2,0.75) node[midway,above,font=\footnotesize]{direction of travel};
\foreach \k in {0,...,6}{\pgfmathtruncatemacro{\j}{\k+1}\node[sl,anchor=north] at (\k*\lw,-0.05) {$b_{v\j}$};}\node[sl,anchor=north] at (7*\lw,-0.05) {$r_{v7}$};
\coordinate (e2) at (2*\lw,0.5); \coordinate (e5) at (5*\lw,0.5);
\fill (e2) circle (2pt); \fill (e5) circle (2pt);
\node[station] (a) at (2*\lw-0.9,2.6) {$c_1$};
\node[station] (b) at (2*\lw+0.9,2.6) {$c_2$};
\node[station] (c) at (5*\lw+0.2,2.6) {$c_3$};
\draw[det] (e2) to[out=110,in=-80] (a.south);
\draw[det,black!45] (e2) to[out=70,in=-100] (b.south);
\draw[det] (e5) to[out=95,in=-95] node[midway,right,font=\footnotesize]{detour and return} (c.south);
\node[font=\footnotesize,align=center] at (2*\lw,3.45) {candidates $\sC_{v2}=\{c_1,c_2\}$};
\end{tikzpicture}
\caption{En-route charging decisions of trip $v$ on a seven-link route. $b_{vi}$ is the SoC at the start of link $\ell_i$, and $r_{v7}$ the SoC on arrival at the destination. At the end of $\ell_2$ the trip may detour to either of two candidate stations; here it uses $c_1$, charges to $\overline{S}$, and returns to the same point. It stops again after $\ell_5$ at its only candidate, $c_3$.}
\label{fig:problem}
\end{figure}

\subsection{Mixed-integer formulation}\label{sec:mip}
\Cref{tab:notation} lists the notation. For trip $v$, link $i$, and candidate $c \in \sC_{vi}$, the binary variable $x_{vic}$ indicates a stop at $c$ after link $i$; $u_{vic} \geq 0$ is the charging duration, $g_{vic} \geq 0$ the waiting time before charging, and $s_{vic}$ the time charging starts. The SoC and the clock at the start of link $i$ are $b_{vi}$ and $\tau_{vi}$. For each station, $z_{vicp}$ assigns a stop to plug $p \in \{1,\ldots,Z_c\}$, and for each pair of stops of different trips at the same station, $y_{vic,v'i'c}$ orders them. Let $r_{vi} = b_{vi} - E_{vi}$ denote the SoC at the end of link $i$ and $\xi_{vi} = \sum_{c \in \sC_{vi}} x_{vic}$. Candidate sets are empty after the last link ($\sC_{vN_v} = \emptyset$), and pairs of stops are ordered by an arbitrary fixed order $\prec$, so each pair has a single sequencing variable; in \crefrange{eq:seq1}{eq:seq2}, $v \neq v'$, $c \in \sC_{vi} \cap \sC_{v'i'}$, and $p \leq Z_c$. Because a stop at which the vehicle would arrive above $\overline{S}$ only adds detour and lowers the SoC, such stops are excluded by the upper bound in~\cref{eq:reach}, which keeps $u_{vic} \geq 0$ in~\cref{eq:dur}.

\begingroup\small
\begin{align}
\min \quad & \sum_{v \in \sV} \sum_{i \in \sL_v} \sum_{c \in \sC_{vi}} \bigl( 2\, C_\tau T^\delta_{vic}\, x_{vic} + C_\chi\, u_{vic} + C_\omega\, g_{vic} \bigr) \label{eq:obj}\\
\text{s.t.}\quad
& b_{v1} = B_v,\quad \tau_{v1} = T^0_v && \forall v \in \sV \label{eq:init}\\
& r_{vi} \geq \underline{S},\quad \xi_{vi} \leq 1 && \forall v \in \sV,\, i \in \sL_v \label{eq:soclow}\\
& \underline{S} - \bigM (1 - x_{vic}) \leq r_{vi} - E^\delta_{vic} \leq \overline{S} + \bigM (1 - x_{vic}) && \forall v \in \sV,\, i \in \sL_v,\, c \in \sC_{vi} \label{eq:reach}\\
& b_{v,i+1} = (1-\xi_{vi})\, r_{vi} + \textstyle\sum_{c \in \sC_{vi}} x_{vic} \bigl(\overline{S} - E^\delta_{vic}\bigr) && \forall v \in \sV,\, i \in \sL_v \setminus \{N_v\} \label{eq:socupd}\\
& u_{vic} = x_{vic}\, \bigl(\overline{S} - r_{vi} + E^\delta_{vic}\bigr) / \phi_v && \forall v \in \sV,\, i \in \sL_v,\, c \in \sC_{vi} \label{eq:dur}\\
& s_{vic} = \tau_{vi} + T_{vi} + T^\delta_{vic} + g_{vic},\quad g_{vic} \leq \bigM x_{vic} && \forall v \in \sV,\, i \in \sL_v,\, c \in \sC_{vi} \label{eq:start}\\
& \tau_{v,i+1} = \tau_{vi} + T_{vi} + \textstyle\sum_{c \in \sC_{vi}} \bigl(2\,T^\delta_{vic}\, x_{vic} + g_{vic} + u_{vic}\bigr) && \forall v \in \sV,\, i \in \sL_v \setminus \{N_v\} \label{eq:clock}\\
& \textstyle\sum_{p=1}^{Z_c} z_{vicp} = x_{vic} && \forall v \in \sV,\, i \in \sL_v,\, c \in \sC_{vi} \label{eq:plug}\\
& s_{vic} + u_{vic} \leq s_{v'i'c} + \bigM\bigl(1-y_{vic,v'i'c}\bigr) + \bigM\bigl(2 - z_{vicp} - z_{v'i'cp}\bigr) && \forall (v,i) \prec (v',i'),\, c,\, p \label{eq:seq1}\\
& s_{v'i'c} + u_{v'i'c} \leq s_{vic} + \bigM\, y_{vic,v'i'c} + \bigM\bigl(2 - z_{vicp} - z_{v'i'cp}\bigr) && \forall (v,i) \prec (v',i'),\, c,\, p \label{eq:seq2}\\
& \underline{S} \leq b_{vi} \leq \max\{B_v, \overline{S}\},\quad x, y, z \in \{0,1\},\quad g, u \geq 0 && \forall v \in \sV,\, i \in \sL_v \label{eq:dom}
\end{align}
\endgroup

The objective~\cref{eq:obj} values the round-trip detour time at $C_\tau$, the charging time at $C_\chi$, and the waiting time at $C_\omega$ per hour; all times are in hours. Constraints~\cref{eq:init} fix the starting state. Constraints~\cref{eq:soclow} keep the SoC above $\underline{S}$ at every link end and allow at most one stop per link endpoint; \cref{eq:reach} requires enough charge to reach the chosen station and excludes stops reached above $\overline{S}$. Constraints~\cref{eq:socupd} set the SoC at the start of the next link to the end-of-link SoC if the trip does not stop and to $\overline{S}$ minus the return-leg energy if it does. Constraints~\cref{eq:dur} set the charging duration to the time needed to reach $\overline{S}$. Constraints~\crefrange{eq:start}{eq:clock} advance each trip's clock through travel, detours, waiting, and charging, so that a delay at one stop shifts all later stops. Constraints~\crefrange{eq:plug}{eq:seq2} assign each stop to one plug and forbid two stops of different trips from overlapping on the same plug. The products in \crefrange{eq:socupd}{eq:dur} are linearized with standard big-$\bigM$ constraints; \cref{sec:benchmark} uses this linear model as the exact benchmark.

\begin{table}[!htb]
\centering
\caption{Notation (alphabetical within each group).}\label{tab:notation}
\footnotesize
\begin{tabularx}{\textwidth}{lX}
\toprule
Symbol & Definition \\
\midrule
\multicolumn{2}{l}{\emph{Sets}}\\
$\sC$ & Charging stations \\
$\sC_{vi}$ & Candidate stations at the end of link $i$ of trip $v$; $|\sC_{vi}| \leq K$ \\
$\sL_v$ & Links of trip $v$, in route order, $\{1,\ldots,N_v\}$ \\
$\sV$ & Electrified trips (vehicles) \\
\midrule
\multicolumn{2}{l}{\emph{Parameters}}\\
$B_v$ & Starting SoC of trip $v$ (\%) \\
$C_\tau$, $C_\chi$, $C_\omega$ & Value of detour time, charging time, and waiting time (\$/h) \\
$\bar{D}$ & Maximum detour distance for candidates other than the nearest station (miles) \\
$E_{vi}$ & SoC used by trip $v$ on link $i$ (\%) \\
$E^\delta_{vic}$ & SoC used by the one-way detour from the end of link $i$ to station $c$ (\%) \\
$K$ & Maximum number of candidate stations per link endpoint \\
$\phi_v$ & Charging rate of trip $v$, $10^5\,W/Q_v$ (\%/h) \\
$Q_v$, $\rho_v$, $R_v$ & Battery capacity (Wh), energy use (Wh/mile), and range (miles) of trip $v$; $Q_v = \rho_v R_v$ \\
$\underline{S}$, $\overline{S}$ & Minimum SoC and target SoC of en-route charging (\%) \\
$\bar{s}$ & Detour speed (mph) \\
$T^0_v$ & Departure time of trip $v$ (h) \\
$T_{vi}$ & Travel time of trip $v$ on link $i$ (h) \\
$T^\delta_{vic}$ & One-way detour time from the end of link $i$ to station $c$ (h) \\
$\vartheta$ & Detour circuity factor applied to $L_1$ distances \\
$W$ & Charger power (kW) \\
$Z_c$ & Number of plugs at station $c$ \\
\midrule
\multicolumn{2}{l}{\emph{Variables}}\\
$b_{vi}$ & SoC of trip $v$ at the start of link $i$ (\%) \\
$g_{vic}$ & Waiting time before charging at station $c$ after link $i$ (h) \\
$s_{vic}$ & Time charging starts at station $c$ after link $i$ (h) \\
$\tau_{vi}$ & Time trip $v$ starts link $i$ (h) \\
$u_{vic}$ & Charging duration at station $c$ after link $i$ (h) \\
$x_{vic}$ & 1 if trip $v$ charges at station $c$ after link $i$ \\
$y_{vic,v'i'c}$ & 1 if stop $(v,i,c)$ precedes stop $(v',i',c)$ on their common plug \\
$z_{vicp}$ & 1 if stop $(v,i,c)$ uses plug $p$ of station $c$ \\
\bottomrule
\end{tabularx}
\end{table}

The MIP cannot be solved at national scale. The national DCFC instance has 2.2 million trips with several hundred links each and up to three candidates per link endpoint, so \crefrange{eq:obj}{eq:dom} has on the order of $10^9$ stop variables before any sequencing variables, and the number of sequencing pairs grows quadratically with the number of stops that can reach the same station. Without the plug-capacity constraints~\crefrange{eq:plug}{eq:seq2}, however, the problem separates by trip: waiting is zero and every trip faces an independent single-vehicle problem. The framework in \cref{sec:method} exploits this structure.
\section{Solution framework}\label{sec:method}

The framework has three stages and a common evaluation (\cref{fig:pipeline}). Preprocessing builds each trip's route, energy use, and candidate stations. A per-vehicle DP then finds each trip's optimal schedule ignoring plug capacity. A congestion stage couples the trips through the plugs. Every schedule, whichever congestion method produced it, is scored with the same FIFO simulation, except where noted (the MIP of \cref{sec:benchmark} is reported by its own objective).

\begin{figure}[!htb]
\centering
\begin{tikzpicture}[box/.style={draw, rounded corners=3pt, minimum height=1.25cm, text width=3.3cm, align=center, font=\footnotesize, fill=black!4},
                    arr/.style={->, >=stealth, thick}]
\node[box] (a) at (0,0) {\textbf{Preprocessing}\\routes, EV attributes, $K$ candidate stations};
\node[box] (b) at (4.1,0) {\textbf{Per-vehicle DP}\\exact shortest path on the charging-opportunity DAG};
\node[box] (c) at (8.2,0) {\textbf{Congestion stage}\\none, station penalties, FIFO best response, or sequential reservation};
\node[box] (d) at (12.3,0) {\textbf{Common evaluation}\\FIFO simulation with delay propagation; economics};
\draw[arr] (a) -- (b); \draw[arr] (b) -- (c); \draw[arr] (c) -- (d);
\end{tikzpicture}
\caption{Solution framework. The DP is exact for the capacity-free problem of each trip; the congestion stage couples trips through plug capacity; all outcomes are scored by the same FIFO simulation.}
\label{fig:pipeline}
\end{figure}

\subsection{Preprocessing}
Each origin--destination trip is routed on the national road network by static shortest-path assignment, which gives its sequence of links with lengths and travel times. Trips with an origin or destination in Alaska or Hawaii are removed. Each trip receives an EV range $R_v \sim \mathrm{Uniform}(\underline{R}, \overline{R})$, an energy-use rate $\rho_v \sim \mathrm{Triangular}(\underline{\rho}, \underline{\rho}, \overline{\rho})$ (most vehicles near the efficient end), and a starting SoC $B_v \sim \mathrm{Triangular}(\underline{B}, \overline{B}, \overline{B})$ (most trips start near full). Link energy is $E_{vi} = 100\, \rho_v D_{vi} / Q_v$ for link length $D_{vi}$ in miles. A ball tree on station coordinates with the $L_1$ metric returns, for every link endpoint on any route, the $K$ nearest stations; detour times and energies follow from (A3). Two filters remove trips that the scheduler need not consider: trips that can finish without charging ($B_v - \sum_{i} E_{vi} \geq \underline{S}$) and trips containing a link that uses more than half of the usable SoC window $\overline{S} - \underline{S}$, a practical filter that removes 15,711 trips (0.6\%). A forward reachability check then removes trips that cannot reach a station of the network before falling below $\underline{S}$ (\cref{sec:design}).

\subsection{Exact per-vehicle dynamic program}\label{sec:dag}
Without plug capacity, trip $v$ solves
\begin{equation}
\min_{x_v} \sum_{i \in \sL_v} \sum_{c \in \sC_{vi}} x_{vic} \Bigl( 2\, C_\tau T^\delta_{vic} + C_\chi \bigl(\overline{S} - r_{vi} + E^\delta_{vic}\bigr)/\phi_v + \pi_c \Bigr) \label{eq:single}
\end{equation}
subject to \crefrange{eq:init}{eq:socupd} and \cref{eq:dom} for that trip, where $\pi_c \geq 0$ is a station penalty that is zero unless the penalty heuristic of \cref{sec:congestion} sets it. Because every stop charges to $\overline{S}$ (A4), the SoC on leaving a stop at $(i, c)$ is $\overline{S} - E^\delta_{vic}$, regardless of what happened before. This observation lets us drop SoC from the state. When $B_v > \overline{S}$, the trip could reach an early station above $\overline{S}$; as in the MIP, such stops are excluded, because skipping one leaves every later stop reachable with at least as much charge and at no greater cost.

Define the \emph{charging-opportunity DAG} $G_v = (\sN_v, \mathcal{A}_v)$ with nodes $\sN_v = \{o\} \cup \{(i, c) : i \in \sL_v, c \in \sC_{vi}\} \cup \{d\}$, where $o$ is the origin and $d$ the destination. An arc from the origin to $(j, c)$ exists if the trip can reach station $c$ after link $j$ on its starting charge without arriving above $\overline{S}$:
\begin{equation}
\underline{S} \leq B_v - \textstyle\sum_{k=1}^{j} E_{vk} - E^\delta_{vjc} \leq \overline{S}. \label{eq:arc0}
\end{equation}
An arc from $(i, c)$ to $(j, c')$ with $j > i$ exists if, after charging at $c$, the trip can drive links $i+1, \ldots, j$ and reach $c'$ without falling below $\underline{S}$:
\begin{equation}
\overline{S} - E^\delta_{vic} - \textstyle\sum_{k=i+1}^{j} E_{vk} - E^\delta_{vjc'} \geq \underline{S}. \label{eq:arc}
\end{equation}
Arcs into the destination are defined in the same way with no final detour, including the arc from $o$ to $d$ for a plan without stops. The cost of an arc into $(j, c')$ is the bracketed term of \cref{eq:single} for that stop, which is known exactly because the arrival SoC $r_{vj} - E^\delta_{vjc'}$ is determined by the arc's tail. Since energy use is nonnegative, the links that can be reached from a node form a contiguous range, which a binary search on cumulative energy finds.

\begin{proposition}\label{prop:exact}
Under \textup{(A3)--(A5)}, every feasible charging plan of trip $v$ in the capacity-free problem corresponds to an origin--destination path in $G_v$ with the same cost, and vice versa. A shortest path in $G_v$ is therefore an optimal solution of the single-trip problem \cref{eq:single}.
\end{proposition}
\begin{proof}
A plan is an ordered set of stops $(i_1, c_1), \ldots, (i_m, c_m)$ with $i_1 < \cdots < i_m$. Between two consecutive stops the trip starts at $\overline{S} - E^\delta_{v i_k c_k}$ and drives links $i_k+1, \ldots, i_{k+1}$, so the plan satisfies the SoC constraints iff each consecutive pair satisfies \cref{eq:arc} (and the first stop satisfies \cref{eq:arc0}, the last the destination condition). Because $E_{vk}, E^\delta_{vic} \geq 0$, SoC does not increase between stops, so the lower bound $\underline{S}$ holds at every link endpoint of a segment if it holds on arrival at the segment's station or destination, which is what the arc conditions check; the return leg of a stop, $\overline{S} - E^\delta_{v i_k c_k} \geq \underline{S}$, is implied by the condition of its outgoing arc. Only a stop reached on the starting charge can arrive above $\overline{S}$, which the upper bound in \cref{eq:arc0} excludes, as the MIP does. The cost of the plan is the sum of the stop costs, each of which depends only on the stop and its predecessor through the arrival SoC; this is the arc cost. Hence plans and paths are in one-to-one correspondence with equal costs. $G_v$ is acyclic because arcs go forward in link index, so a forward pass in topological order solves it.
\end{proof}

With $m_v = \sum_i |\sC_{vi}|$ opportunities, the forward pass relaxes, for each node, the opportunities reachable within one charge, so its cost is $O(m_v \cdot \bar{a}_v)$, where $\bar{a}_v$ is the average number of opportunities within one charge. In contrast, a DP over (link, SoC bucket) states has $O(N_v\, K\, (\overline{S}-\underline{S})/\gamma)$ states and must round each link's energy to a whole number of buckets of size $\gamma$. Rounding each link up to at least one bucket overstates energy on short links. National routes have hundreds of links, so this overstatement is large (\cref{sec:ablation}). The opportunity DAG avoids it entirely. \Cref{alg:dag} summarizes the DP; the trips are solved independently in parallel.

\begin{algorithm}[!htb]
\caption{Exact per-vehicle DP on the charging-opportunity DAG}\label{alg:dag}
\SetArgSty{textnormal}
\footnotesize
\KwIn{Route links $\sL_v$ with energies $E_{vi}$; candidates $\sC_{vi}$ with $T^\delta_{vic}$, $E^\delta_{vic}$; $B_v$, $\underline{S}$, $\overline{S}$, $\phi_v$; costs $C_\tau$, $C_\chi$; penalties $\pi$}
\KwOut{Optimal stops $\{(i_k, c_k)\}$ and cost, or \emph{infeasible}}
List opportunities $q = (i_q, c_q)$ in order of $i_q$; compute cumulative energy $\Lambda_j = \sum_{k \leq j} E_{vk}$\;
$\mathrm{dist}(q) \leftarrow \infty$ for all $q$; relax from the origin with SoC $B_v$ at link 1\;
\If{$B_v - \Lambda_{N_v} \geq \underline{S}$}{best $\leftarrow 0$ (no stop needed)}
\ForEach{opportunity $q$ in order with $\mathrm{dist}(q) < \infty$}{
  $\beta \leftarrow \overline{S} - E^\delta_{v i_q c_q}$ \tcp*{SoC when leaving the stop}
  \If{$\beta - (\Lambda_{N_v} - \Lambda_{i_q}) \geq \underline{S}$}{update best with $\mathrm{dist}(q)$, predecessor $q$}
  Find by binary search the last link $j$ with $\beta - (\Lambda_j - \Lambda_{i_q}) \geq \underline{S}$\;
  \ForEach{opportunity $w$ with $i_q < i_w \leq j$ and $\beta - (\Lambda_{i_w} - \Lambda_{i_q}) - E^\delta_{v i_w c_w} \geq \underline{S}$}{
    $a \leftarrow 2C_\tau T^\delta_{v i_w c_w} + \pi_{c_w} + C_\chi\,(\overline{S} - \beta + \Lambda_{i_w} - \Lambda_{i_q} + E^\delta_{v i_w c_w})/\phi_v$\;
    \lIf{$\mathrm{dist}(q) + a < \mathrm{dist}(w)$}{$\mathrm{dist}(w) \leftarrow \mathrm{dist}(q) + a$; $\mathrm{pred}(w) \leftarrow q$}
  }
}
\Return stops from predecessors of the best terminal, or \emph{infeasible} if none\;
\end{algorithm}

\subsection{Congestion stage}\label{sec:congestion}
The per-vehicle schedules ignore other trips. We study four ways of coupling them: no coordination, a station-penalty heuristic, FIFO best response, and sequential reservation.

Under no coordination (DAG + FIFO), each trip keeps its capacity-free optimal schedule and queues FIFO wherever plugs are busy; in the national case runs a two-hour best-response stage follows, which changed less than 0.5\% of the cost. This is the outcome if travelers plan in isolation, and it gives the reference against which coordination is measured.

The station-penalty heuristic adjusts station costs to the observed overflow. The plan's sessions are binned into 15-minute intervals; a station whose occupancy exceeds $Z_c$ in some bin receives a penalty $\pi_c = \lambda \max_t (\text{occupancy}_{ct} - Z_c)^+$; the trips that use penalized stations are re-solved with \cref{alg:dag}; and the steps repeat, with $\lambda$ doubled when overflow stops decreasing, until no overflow remains, an iteration limit is reached, or a time limit expires. The final iterate is then queued FIFO. Penalties are recomputed from the current overflow in each iteration rather than accumulated, so a station that is vacated in one iteration loses its penalty in the next. The heuristic has no convergence guarantee and is not a dual method.

FIFO best response starts from the capacity-free schedules and simulates all trips FIFO with delay propagation. Each trip that waited is then removed in turn, in decreasing order of waiting time, and re-planned against the realized plug occupancy of all other trips by a time-dependent version of \cref{alg:dag}: when evaluating a stop, the DP reads the earliest time a plug at that station is free for the needed duration from a per-plug reservation calendar, adds the waiting time to the arc cost at $C_\omega$, and carries the resulting departure time forward. The trip's new plan is accepted if it lowers its own cost, and its sessions are re-inserted. After each round the system is re-simulated and the best system state is kept, so the result is never worse than DAG + FIFO. In the time-dependent DP each node keeps only its least-cost label, so the DP is exact without waiting and heuristic with it.

Sequential reservation plans trips one at a time in order of departure. Each trip solves the time-dependent DP of best response against a per-plug reservation calendar that holds the sessions of all earlier trips, and its own sessions are then added to the calendar. Waiting enters the trip's objective at $C_\omega$ when it plans, so a trip whose preferred station is booked stops at another candidate station, at another link endpoint, or accepts a short wait, whichever costs least. If a trip's capacity-free schedule does not wait under the current calendar, it is accepted without solving the time-dependent DP; on all trips at 350~kW this fast path covers about 44\% of them. Back-to-back sessions on a plug are stored as merged blocks, so the earliest free time at a busy station is found in time proportional to the number of gaps rather than the number of sessions; this changes running time, not results. The scheme works like a first-come, first-served reservation service: in the plan, a trip is never displaced by a trip that departs later (the common evaluation can still record waits; \cref{sec:S3}), so the outcome depends on the order and is not a system optimum. A best-response round can follow the pass. The difference from FIFO best response is the starting point. Best response starts from a full FIFO simulation in which millions of trips already wait and must undo that congestion one trip at a time. The reservation pass never creates the congestion, so each trip is planned once.

Independently of the congestion method, the candidate sets can also account for station size. The $K$ nearest stations of a link endpoint can all be one- and two-plug sites even when a large station is a short detour away. The capacity-aware variant keeps the nearest station and fills the other $K-1$ slots with the stations that have the most plugs among the 20 nearest within $\bar{D}$, breaking ties by distance. It changes only which stations a trip may consider; the DP and the congestion stage are unchanged.

\subsection{Evaluation of schedules and costs}\label{sec:eval}
To compare congestion methods on equal terms, every final plan is scored, except where noted, with one event simulation. Trips depart at $T^0_v$ and drive to their first stop. Stops are served in order of arrival at each station on the plug that frees first. A trip that waits starts charging late, and all of its later arrivals are shifted by the delay. The MIP is less restrictive: it may order stops on a plug in any sequence and hold a plug for a later arrival, so its solutions are reservation schedules, and its objective can be lower than the FIFO cost of the same stops. In \cref{sec:benchmark} we report the MIP objective as is, which favors the MIP. The penalty heuristic's own accounting queues each stop at its planned, undelayed arrival time and is used only to describe that heuristic.

\label{sec:econ}%
For each trip we compare the EV with an ICE vehicle driving the same route, at three cost boundaries. Let $d_v$ be the route length and $d^\delta_v$ the total detour distance (miles), and let $t^\delta_v$, $t^\chi_v$, and $t^\omega_v$ be the trip's total detour, charging, and waiting times from this evaluation.
\begin{align}
\text{Tier 1 (operating cost):}\quad & C^{\mathrm{ICE},1}_v = \frac{d_v}{\eta} P_g + C_\chi \Bigl\lfloor \frac{d_v}{d_g} \Bigr\rfloor t_g, \notag\\
& C^{\mathrm{EV},1}_v = \frac{(d_v + d^\delta_v)\rho_v}{1000} P_e + C_\tau t^\delta_v + C_\chi t^\chi_v + C_\omega t^\omega_v, \label{eq:tier1}\\
\text{Tier 2 (ownership):}\quad & C^{\mathrm{ICE},2}_v = C^{\mathrm{ICE},1}_v + \Bigl(\frac{V}{L\,M} + m_g\Bigr) d_v, \notag\\
& C^{\mathrm{EV},2}_v = C^{\mathrm{EV},1}_v + \Bigl(\frac{\alpha V}{L\,M} + m_e\Bigr) (d_v + d^\delta_v), \label{eq:tier2}\\
\text{Tier 3 (infrastructure):}\quad & C^{\mathrm{EV},3}_v = C^{\mathrm{EV},2}_v + \frac{4.345\, \sigma\, \kappa\, Z\, (1/L_c + \mu)}{52\, |\sV|}, \label{eq:tier3}
\end{align}
where $\eta$ is ICE fuel economy, $P_g$ and $P_e$ are gasoline and electricity prices, $t_g$ is the time of one ICE refueling stop every $d_g$ miles, $V$ is the ICE purchase price, $\alpha$ the EV price ratio, $L$ the vehicle life in years and $M$ the annual mileage (so vehicle cost is amortized per mile), $m_g$ and $m_e$ are maintenance costs per mile, $\kappa$ is the cost per plug, $Z$ the number of plugs, $L_c$ the charger life, $\mu$ the annual operations-and-maintenance share, and $\sigma$ the sample fraction; the annual plug cost is spread over a year of national long-distance trips, $|\sV|/\sigma$ per month (52 weeks per year, 4.345 per month). We report the ratio-of-totals saving $100\,(\sum_v C^{\mathrm{ICE}}_v - \sum_v C^{\mathrm{EV}}_v)/\sum_v C^{\mathrm{ICE}}_v$, which is not dominated by short trips as a mean of per-trip percentages would be, without the waiting term $C_\omega t^\omega_v$ and, where waits are bounded, with it.
\section{Data and experimental design}\label{sec:design}

On the supply side, the Alternative Fuels Station Locator \citep{afdc_2026} lists 14,260 public direct-current fast-charging (DCFC) stations with 68,362 plugs and 64,477 public Level~2 (L2) stations with 169,385 plugs in the contiguous United States. The two networks are studied separately. Plug counts are taken as reported; 1- and 2-plug stations are common in both networks. More than half of the DCFC plugs (54\%) belong to the Tesla network. Because every plug is treated as usable by every EV and every station as delivering the scenario's single power level, the 50 and 350~kW scenarios bound the charging speed of the mixed network rather than describe it.

On the demand side, POLARIS synthesizes one month of national long-distance personal-vehicle trips for a 1\% sample of the U.S.\ population, 2,768,665 trips after removing Alaska and Hawaii \citep{auld_2016}. POLARIS applies no fixed distance threshold to long-distance trips; 98\% of the automobile trips are 50 miles or longer, the threshold of the National Household Travel Survey. Departure times follow the simulated activity schedules and span the whole month. Every trip is assumed to be made by an EV. The simulated automobile trips are longer than survey data suggest: their median is 392 miles and 37\% are 500 miles or longer, whereas in the 2001 National Household Travel Survey about 10\% of long-distance trips by all modes reached 500 miles \citep{bts_2006}; that figure includes air travel, so the share among automobile trips is likely lower still. The sample is therefore weighted toward the trips that need the most en-route charging, which makes charging loads high and the cost comparison unfavorable to the EV. After the preprocessing filters, which remove the 15,711 trips with an overly long link and the trips that need no charge, 2,244,001 trips cannot finish on their starting charge; of these, the forward check screens out 33,428 (1.5\%) that do not reach a DCFC station along their route before falling below $\underline{S}$ and 74,002 (3.3\%) that do not reach an L2 station; L2 stations are more numerous but, presumably because they cluster in urban areas, leave longer gaps along some rural highways. The check stops at the last station within reach, so it may also count a few trips that an earlier stop would save. These trips are excluded before scheduling, which leaves 2,210,573 trips on the DCFC network and 2,169,999 on the L2 network. The counts therefore differ between networks but not between power levels, which do not affect reachability; the feasible counts of \cref{tab:case} can still differ by a trip, through the few edge cases discussed there. The full data set therefore represents electrifying every long-distance trip of a 1\% sample, or about 1\% of national long-distance demand, served by the entire national station inventory. Neither demand nor supply is rescaled: the experiment asks how the existing network handles a demand level of the same order as today's EV share of the light-duty fleet if those trips were scheduled as modeled, and \cref{sec:discussion} discusses what full adoption would imply.

The case study crosses the two networks with two power levels each: DCFC at 50~kW (the lowest DCFC power) and 350~kW, and L2 at 5 and 15~kW. The L2 scenarios are not a proposal to serve long trips with L2; they bound what the much denser L2 network could contribute if it were the only option. All four use the parameters of \cref{tab:params}, three candidate stations per link endpoint, and $\bar{D} = 15$ miles. Sensitivity sweeps use a random 1\% subsample of the trips (27,687 trips, of which about 22,000 need to charge at baseline values; about 0.01\% of national demand). Each level is replicated with ten random seeds that redraw the EV attributes.

\begin{table}[!htb]
\centering
\caption{Parameters of the experiments.}\label{tab:params}
\footnotesize
\begin{tabularx}{\textwidth}{l X l >{\raggedright\arraybackslash}p{3.6cm}}
\toprule
Symbol & Meaning (unit) & Baseline & Source \\
\midrule
\multicolumn{4}{l}{\emph{Vehicles and charging}}\\
$\underline{R}$, $\overline{R}$ & EV range bounds, uniform (miles) & 150, 350 & \citealp{iea_2024} \\
$\underline{\rho}$, $\overline{\rho}$ & Energy-use bounds, triangular with mode $\underline{\rho}$ (Wh/mile) & 250, 750 & \citealp{iea_2024} \\
$\underline{B}$, $\overline{B}$ & Starting-SoC bounds, triangular with mode $\overline{B}$ (\%) & 60, 100 & Assumed \\
$\underline{S}$, $\overline{S}$ & Minimum SoC, target SoC of en-route charging (\%) & 20, 80 & Assumed \\
$W$ & Charger power (kW) & 50, 350 (DCFC); 5, 15 (L2) & \citealp{afdc_2026} \\
$K$, $\bar{D}$ & Candidate stations per link endpoint; distance limit for non-nearest candidates (miles) & 3, 15 & Assumed \\
$\vartheta$, $\bar{s}$ & Detour circuity factor; detour speed (mph) & 1.0, 40 & Assumed \\
\midrule
\multicolumn{4}{l}{\emph{Values of time}}\\
$C_\tau$ & Detour time (\$/h) & 70 & Assumed \\
$C_\chi$ & Charging time and ICE refueling time (\$/h) & 50 & Assumed \\
$C_\omega$ & Waiting time (\$/h) & 50 & Assumed \\
\midrule
\multicolumn{4}{l}{\emph{Economics}}\\
$\eta$, $P_g$ & ICE fuel economy (mpg); gasoline price (\$/gal) & 25; 3.50 & \citealp{epa_trends_2025,eia_gas_2026} \\
$d_g$, $t_g$ & ICE range between refueling stops (miles); refueling stop (min) & 350; 10 & Assumed \\
$P_e$ & Electricity price (\$/kWh); en-route sensitivity & 0.15; 0.25--0.45 & \citealp{eia_elec_2024}; assumed \\
$V$, $\alpha$ & ICE purchase price (\$); EV/ICE price ratio & 35,000; 1.2 & \citealp{iea_2024} \\
$L$, $M$ & Vehicle life (years); annual mileage (miles) & 12; 13,500 & Assumed \\
$m_g$, $m_e$ & Maintenance cost, ICE and EV (\$/mile) & 0.09, 0.04 & \citealp{aaa_2024} \\
$\kappa$, $L_c$, $\mu$ & Cost per plug (\$); charger life (years); annual O\&M (share of capital) & 50,000; 12; 0.08 & Assumed \\
$\sigma$ & Demand sample fraction & 0.01 & POLARIS sample \\
\bottomrule
\end{tabularx}
\par\smallskip
\begin{minipage}{\linewidth}\footnotesize\textit{Note:} Baseline values are used unless a sweep varies them.\end{minipage}
\end{table}

The values of time in \cref{tab:params} are deliberately round and per vehicle-hour. The U.S.\ DOT guidance for personal intercity travel, \$19.00 per person-hour in 2015 dollars~\citep{usdot_vot_2024}, is about \$54 per vehicle-hour in 2026 dollars for two occupants; our charging and waiting values are close to it and our detour value is above it. Time spent detouring, charging, and waiting at a station is treated as largely unproductive. Because the conclusions depend on the value of charging time, \cref{sec:econresults} reports the break-even value at which EV and ICE costs are equal, so readers can apply their own.

The validation in \cref{sec:validation} uses the DCFC 50~kW scenario on the full data set and, separately, small synthetic corridor instances described in \cref{sec:benchmark}; \cref{sec:national} reports the four national scenarios.

\section{Validation}\label{sec:validation}

\subsection{Comparison with the exact MIP on corridor instances}\label{sec:benchmark}
We generated 25 corridor instances small enough for the exact model~\crefrange{eq:obj}{eq:dom}: 5 to 30 trips that depart within the same hour share a 600-mile highway with 18 single-plug stations (\hyperref[sec:supp]{Supplementary material}, \cref{sec:S1}). These instances are deliberately congested, a harder setting for coordination than the national data, whose stations have more plugs and whose departures spread over a month. The MIP is solved with HiGHS (one thread, 120~s limit), and all heuristic schedules are scored with the evaluation of \cref{sec:eval}, so all heuristic costs are on the same basis; the MIP is reported as its own objective (\cref{sec:eval}). We compare the four congestion methods of \cref{sec:congestion}, all with the nearest candidate sets. \Cref{tab:bench} reports the average gap to the best solution found by any method.

\begin{table}[!htb]
\centering
\begin{threeparttable}
\caption{Average gap (\%) to the best known solution on 25 congested corridor instances.}\label{tab:bench}
\footnotesize
\setlength{\tabcolsep}{4pt}\begin{tabular}{rcrrrrrrrr}
\toprule
Trips & Proven & MIP (s) & DAG+FIFO & Penalty & Seq.\ res. & BR & MIP inc. & BR (s) & BR vs.\ LB \\
\midrule
5  & 5/5 & 1   & 3.9  & 13.6 & 10.3 & 2.2 & 0.0  & 0.01 & 2.2 \\
10 & 4/5 & 68  & 14.5 & 34.5 & 4.7  & 5.6 & 0.0  & 0.04 & 5.9 \\
15 & 2/5 & 78  & 15.4 & 39.4 & 4.6  & 1.8 & 7.6  & 0.06 & 6.9 \\
20 & 0/5 & 120 & 36.2 & 74.6 & 3.3  & 1.1 & 26.7 & 0.10 & 20.7 \\
30 & 0/5 & 120 & 28.1 & 65.6 & 7.4  & 0.2 & 39.6 & 0.13 & 98.5 \\
\midrule
All & 11/25 & 78 & 19.6 & 45.5 & 6.1 & 2.2 & 14.8 & 0.07 & 26.9 \\
\bottomrule
\end{tabular}
\begin{tablenotes}[flushleft]\footnotesize
\item \textit{Note:} Five instances per size. Proven: instances the MIP solved to optimality within 120~s. Heuristic schedules are scored with the same FIFO simulation with delay propagation; the MIP incumbent is its objective value (a reservation schedule), so the best known solution may be a value that is not scored by the FIFO simulation, which makes the heuristic gaps conservative. Seq.\ res.: sequential reservation; BR: FIFO best response; MIP inc.: MIP incumbent. BR vs.\ LB: gap between FIFO best response and the MIP lower bound, an upper bound on its optimality gap.
\end{tablenotes}
\end{threeparttable}
\end{table}

Four results stand out. First, FIFO best response is within 2.2\% of the best known solution on average (median 0\%, maximum 11.2\%), never worse than DAG + FIFO by construction, and better than the MIP's 120-second incumbent on 11 of 25 instances, in a tenth of a second. The best known solution is often best response itself, and the MIP lower bounds are weak for 20 or more trips, so its distance from optimality is proven small only on the smaller instances (2--7\% on average for up to 15 trips). Second, the uncoordinated outcome is 19.6\% above the best on average and up to 53\%, so on congested corridors coordination is worth a lot. Third, the station-penalty heuristic is worse than doing nothing: it raises cost above DAG + FIFO on 21 of 25 instances and is on average 45.5\% above the best (maximum 162\%). The reason is visible in its mechanics. Penalties are recomputed from the current overflow, so trips pushed off a station in one iteration return when its penalty disappears in the next; and the overflow it reduces counts vehicles per 15-minute bin, which is not the waiting that the vehicles actually experience once delays propagate. Moving trips to farther stations adds detour time without reliably reducing that waiting. Fourth, sequential reservation is 6.1\% above the best on average (median 5.3\%, maximum 18.0\%) and better than no coordination on 19 of 25 instances. On average it is worse than best response at every size but ten trips and worse than no coordination on the five-trip instances, where little waiting arises and its fixed planning order can commit early trips to stations that later trips need more. On the larger instances it is far better than no coordination (3.3--7.4\% against 15--36\%). Its value is that it scales: \cref{sec:coordination} applies it to all national trips at 350~kW, where best response could not.

\subsection{Stepwise comparison of modeling choices}\label{sec:ablation}
To see how much the modeling choices of this paper matter at national scale, we ran the DCFC 50~kW scenario on the full data set five times, starting from a simplified baseline that we built from choices common in large-scale models and adding our choices one at a time. Run A is that baseline: the SoC grid ends at the charging target $\overline{S} = 80\%$, so a trip that starts above it is rejected; SoC is discretized into 1\% buckets; charging time is not valued; and congestion is handled by the station-penalty stage. Run B extends the grid to the starting SoC, so each trip starts at its own SoC $B_v$, drawn between 60\% and 100\% (\cref{tab:params}), as (A4) allows. Run C replaces the SoC-bucket DP with the exact opportunity DAG of \cref{sec:dag}. The bucket DP carries SoC as a state, rounding the energy of every link up to a whole 1\% bucket; the DAG needs no SoC state, because every stop charges to $\overline{S}$ and the SoC after each stop is therefore known exactly. Run D adds the value of charging time $C_\chi$ to the DP and to the economics. Run E replaces station penalties with FIFO best response; it differs from the DCFC 50~kW case study only in using one candidate station per link endpoint instead of three. Runs A--E all use $K = 1$. \Cref{tab:ablation} reports the results. To put waiting on a common basis, we recomputed it for runs B--E from each saved schedule with an independent implementation of the FIFO simulation with delay propagation; for run E and the four case-study scenarios, which already use that simulation, the recomputed totals match the scheduler's own to within rounding, which checks the implementation.

\begin{table}[!htb]
\centering
\caption{Stepwise comparison of modeling choices on the national DCFC 50~kW data set.}\label{tab:ablation}
\footnotesize
\setlength{\tabcolsep}{4pt}
\resizebox{\textwidth}{!}{\begin{tabular}{l l r r r r r r r}
\toprule
Run & Change & Feasible & Stops & Detour & Charge & Tier 1, no & Stations & Trips \\
 & & trips & per trip & (min) & (min/stop) & wait (\%) & over capacity & wait $>$1 h (\%) \\
\midrule
A & Simplified baseline & 644,279 & 8.47 & 38.6 & 29.0 & $-$10.8 & --- & --- \\
B & + start SoC above $\overline{S}$ & 2,209,450 & 8.67 & 56.5 & 28.4 & $-$37.5 & 684 & 80.1 \\
C & + exact opportunity DAG & 2,210,571 & 3.76 & 17.2 & 54.9 & +27.1 & 493 & 70.0 \\
D & + value of charging time & 2,210,570 & 3.72 & 21.0 & 52.3 & $-$170.0 & 398 & 61.0 \\
E & $-$ penalties, + best response & 2,210,570 & 3.65 & 7.7 & 51.6 & $-$143.6 & 682 & 72.8 \\
\bottomrule
\end{tabular}}
\par\smallskip
\begin{minipage}{\linewidth}\footnotesize\textit{Note:} 2,210,573 trips; $K = 1$. Each run adds one modeling choice to the previous row. Tier~1: EV-versus-ICE operating-cost reduction excluding waiting (for A--C, also excluding charging and refueling time). Last two columns: FIFO simulation of \cref{sec:eval} applied to each saved schedule (not retained for run A).\end{minipage}
\end{table}

Three of the four changes alter the national picture by a large factor. Rejecting the trips that start above $\overline{S}$ (run A) makes 71\% of the trips infeasible, because about three-quarters of all trips start above 80\% under the starting-SoC distribution of \cref{tab:params}, and trips that start below it are more likely to need charging; a model that drops such trips without reporting them loses them without notice, whereas we report the trips our reachability check excludes (\cref{sec:design}). Extending the grid (B) makes all but about 1,100 trips feasible. Replacing the SoC buckets with the exact opportunity DAG (C) then cuts the number of stops from 8.7 to 3.8 per trip. The bucket DP rounds the energy of every link up to a whole 1\% bucket. A mile of travel uses only 0.3--0.7\% of the battery of the vehicles simulated here, and trajectories consist of hundreds of short links, so the rounding greatly overstates the energy the DP believes a trip needs and forces stops that are not needed; the energy actually consumed, computed from distances, changes by only 5\%. The error comes from rounding each link separately; finer buckets or rounding the cumulative energy would reduce it at a cost in states, whereas the opportunity DAG needs neither. Excluding charging and refueling time, as runs A--C do, the EV then goes from 37.5\% more expensive than the ICE vehicle to 27.1\% cheaper, mainly through the shorter detours of fewer stops, part of which is penalty-stage detour; the cleaner evidence of the bucket error is the number of stops. Valuing charging time (D), which also adds that time to the cost, reverses it again, to 170\% more expensive, because at 50~kW the 3.2 hours of charging per trip in run D cost far more than the fuel saved; unlike the previous step, this reversal comes from a single modeling choice.

The last change replaces station penalties with best response. Neither method coordinated reliably at this scale: best response re-planned only 560 trips in its two hours, and the penalty stage stopped after a single iteration that re-planned about two million trips at once and added about 13 minutes of detour per trip (\hyperref[sec:supp]{Supplementary material}, \cref{sec:S2}). Run E is therefore effectively the capacity-free reference; the detours of runs B--D include this penalty-stage detour, so the detour changes from B to D are partly confounded with it.
\section{National results}\label{sec:national}

\subsection{Feasibility, stops, and detours}\label{sec:feas}
\Cref{tab:case} summarizes the four scenarios. Of the trips scheduled (\cref{sec:design}), all but two to four per scenario can be completed; the exceptions are edge cases that pass the forward check of \cref{sec:design} but for which the DP finds no feasible plan. Trips that need charging are about 525 miles long on average and make 3.6--3.7 stops. Detours are short: 6.6--7.7 minutes and about 1\% of added distance per trip, because the nearest station is usually close to the highway. The two networks need about the same number of stops because the stopping pattern is set by range and SoC limits, not by station density, but about twice as many L2 as DCFC stations are used, so DCFC trips are concentrated on fewer stations. In capacity-free plans the number of candidates per stop hardly matters at this scale: run E of the stepwise comparison, with one candidate, agrees with the 50~kW case study, with three, to the reported precision on feasibility, stops, detour, charging time, and cost. Outside dense urban areas the second and third nearest stations are rarely worth the extra detour. What differs by orders of magnitude between scenarios is the time spent charging: about 52 minutes per stop at 50~kW and 7 minutes at 350~kW, against 2.8 hours at 15~kW and 8.5 hours at 5~kW L2.

\begin{table}[!htb]
\centering
\begin{threeparttable}
\caption{Trips served through the existing networks (DAG + FIFO).}\label{tab:case}
\footnotesize
\begin{tabular}{l r r r r}
\toprule
 & DCFC 50 kW & DCFC 350 kW & L2 5 kW & L2 15 kW \\
\midrule
Trips scheduled & 2,210,573 & 2,210,573 & 2,169,999 & 2,169,999 \\
Feasible trips & 2,210,570 & 2,210,571 & 2,169,995 & 2,169,995 \\
Stations used & 9,684 & 8,764 & 18,451 & 17,931 \\
Stops per trip & 3.65 & 3.69 & 3.56 & 3.56 \\
Detour (min/trip) & 7.7 & 6.6 & 7.7 & 6.9 \\
Added distance (\%) & 0.98 & 0.84 & 1.00 & 0.89 \\
Charging time (min/stop) & 51.6 & 7.4 & 510.2 & 170.1 \\
Charging time (h/trip) & 3.14 & 0.46 & 30.23 & 10.09 \\
\bottomrule
\end{tabular}
\begin{tablenotes}[flushleft]\footnotesize
\item \textit{Note:} Values per trip are means over trips. DAG + FIFO follows a two-hour best-response stage that changed less than 0.5\% of the cost.
\end{tablenotes}
\end{threeparttable}
\end{table}

\subsection{Travel cost relative to ICE vehicles}\label{sec:econresults}
\Cref{tab:econ} compares the cost of the EV trips with the same trips made by an ICE vehicle, excluding waiting time, which we treat separately in \cref{sec:congestionresults}. ICE costs differ slightly between networks because their trip sets differ (\cref{sec:design}). At the baseline electricity price of \$0.15 per kWh, which is close to residential rates and therefore a lower bound on what public charging costs, energy favors the EV: electricity for the average trip costs \$33, against \$74 of gasoline. Detours add \$7.7--9.0. The cost that decides the comparison is charging time. At 50~kW the average trip spends 3.1 hours charging, worth \$157 at $C_\chi = \$50$ per hour, against 10 minutes of refueling for the ICE vehicle; the EV trip costs 144\% more than the ICE trip. At 350~kW charging takes 28 minutes per trip, and the EV is 22\% cheaper in total; with the coordinated plans of \cref{sec:coordination}, in which waiting is negligible, it is cheaper for 84\% of the individual trips. At L2 power the EV trip costs 6.8--19.4 times as much as the ICE trip (reductions below $-100\%$ mean the EV costs more than twice as much), so L2 stations are not a substitute for DCFC on long trips, whatever their coverage.

Because the comparison hinges on $C_\chi$, \cref{tab:econ} also reports the break-even value of charging time: the value a traveler places on an hour spent charging (not a price paid for charging) at which the EV and ICE trips cost the same in total. At \$0.15 per kWh it is about \$11 per hour at 50~kW and \$111 per hour at 350~kW, so the EV trip is cheaper for a traveler who values time below these levels; at the DOT value for personal intercity travel (\$19 per person-hour in 2015 dollars, about \$27 in 2026 dollars, or \$54 per vehicle-hour for two occupants), it is cheaper at 350~kW and more expensive at 50~kW. If detour time is valued at the same rate as charging time, the Tier~1 reduction becomes $-24$\% at 50~kW and $+43$\% at 350~kW at \$19 per hour, and $-55$\% and $+38$\% at \$27 per hour; at about \$54 per vehicle-hour (two occupants) it becomes $-153$\% and $+23$\%. At 350~kW the break-even value is high because charging takes little time. The comparison is just as sensitive to the price of the energy bought en route. About 157--160 kWh of the 221 kWh an average DCFC trip consumes is charged en route; the rest comes from the starting charge, which we keep at \$0.15. If en-route energy costs \$0.25 per kWh, the break-even value falls to \$5 and \$57 per hour; at \$0.35 the EV trip is cheaper only if charging time is worth almost nothing; and at \$0.45, a common price at public fast chargers, the EV trip costs more than the ICE trip at either power even if charging time is free. On long trips that charge mostly en route, the EV's operating-cost advantage therefore depends on fast-charging prices close to residential rates as much as on charging speed. Adding vehicle ownership and maintenance (Tier~2) adds about the same amount to both vehicles, because the EV's higher purchase price is offset by its lower maintenance cost; it shrinks the percentage differences by enlarging the base but changes neither their sign nor the DCFC break-even values by more than about 7\%. Tier~3 spreads the annualized cost of all plugs over the national long-distance trips that charge en route and changes the results by at most 0.3 percentage points; if instead the cost of the whole network were charged to the 1\% of demand simulated here, it would add about \$21 per DCFC trip. Meal and rest breaks can absorb part of the charging time. If every stop of both vehicles includes a 30-minute break that the traveler would take anyway, the breaks absorb about 57\% of the charging time at 50~kW and all of it at 350~kW, and the Tier~1 reduction becomes about $-50$\% at 50~kW and $+45$\% at 350~kW. Charging every final stop to $\overline{S}$ overstates charging slightly: trips arrive with 24\% (50~kW) and 28\% (350~kW) of charge on average, against $\underline{S} = 20\%$. If the last stop charged only what the trip needs, charging time would fall by about 3\% and 5\% (energy consumed, which \cref{eq:tier1} prices at \$0.15 per kWh, is unchanged), the Tier~1 reduction would become $-139$\% and $+23$\%, and the break-even values about \$11 and \$120 per hour. The long trips of the sample tilt the comparison against the EV: at 350~kW the EV trip is 28\% cheaper for trips under 200 miles, 27\% for trips of 200--500 miles, and 20\% for trips of 500 miles or more, which make up 43\% of the trips that charge; at 50~kW the three values are $-50$\%, $-106$\%, and $-164$\%.

\begin{table}[!htb]
\centering
\begin{threeparttable}
\caption{EV-versus-ICE comparison per trip, excluding waiting.}\label{tab:econ}
\footnotesize
\begin{tabular}{l r r r r}
\toprule
 & DCFC 50 kW & DCFC 350 kW & L2 5 kW & L2 15 kW \\
\midrule
ICE fuel (\$/trip) & 73.5 & 73.5 & 72.2 & 72.2 \\
ICE refueling time (\$/trip) & 8.1 & 8.1 & 7.9 & 7.9 \\
EV electricity (\$/trip) & 33.1 & 33.1 & 32.5 & 32.5 \\
EV detour time (\$/trip) & 9.0 & 7.7 & 9.0 & 8.0 \\
EV charging time (\$/trip) & 156.8 & 22.9 & 1,511.4 & 504.4 \\
\midrule
Tier 1 reduction (\%) & $-$143.6 & +22.0 & $-$1,839.1 & $-$580.4 \\
Tier 2 reduction (\%) & $-$47.5 & +8.3 & $-$618.4 & $-$194.5 \\
Tier 3 reduction (\%) & $-$47.6 & +8.2 & $-$618.6 & $-$194.8 \\
Break-even $C_\chi$, Tier 1 (\$/h) & 10.6 & 110.8 & 1.0 & 3.2 \\
Break-even $C_\chi$, Tier 2 (\$/h) & 11.2 & 118.4 & 1.1 & 3.4 \\
\midrule
\multicolumn{5}{l}{\emph{En-route energy at public charging prices}}\\
En-route energy charged (kWh/trip) & 157 & 160 & 151 & 151 \\
Energy consumed (kWh/trip) & 221 & 221 & 217 & 217 \\
Tier 1 reduction, \$0.25/kWh (\%) & $-$162.8 & +2.4 & $-$1,858.0 & $-$599.3 \\
Break-even $C_\chi$, \$0.25/kWh (\$/h) & 5.3 & 56.5 & 0.5 & 1.7 \\
Tier 1 reduction, \$0.35/kWh (\%) & $-$182.0 & $-$17.3 & $-$1,876.8 & $-$618.2 \\
Break-even $C_\chi$, \$0.35/kWh (\$/h) & 0.0 & 2.2 & 0.0 & 0.1 \\
Tier 1 reduction, \$0.45/kWh (\%) & $-$201.2 & $-$36.9 & $-$1,895.7 & $-$637.1 \\
Break-even $C_\chi$, \$0.45/kWh (\$/h) & none & none & none & none \\
\bottomrule
\end{tabular}
\begin{tablenotes}[flushleft]\footnotesize
\item \textit{Note:} Reductions are $(\text{ICE} - \text{EV}) / \text{ICE}$ of the summed costs; negative values mean the EV is more expensive. Break-even $C_\chi$: value of the traveler's charging time (\$ per hour, not a charging price) at which EV and ICE costs are equal (``none'': the EV costs more even if charging time is free). Upper block: all electricity at \$0.15/kWh. Lower block: energy charged en route at the stated public price, the starting charge at \$0.15/kWh.
\end{tablenotes}
\end{threeparttable}
\end{table}

\subsection{Congestion without coordination}\label{sec:congestionresults}
\begin{figure}[!htb]
\centering
\begin{subfigure}{0.49\textwidth}\includegraphics[width=\linewidth]{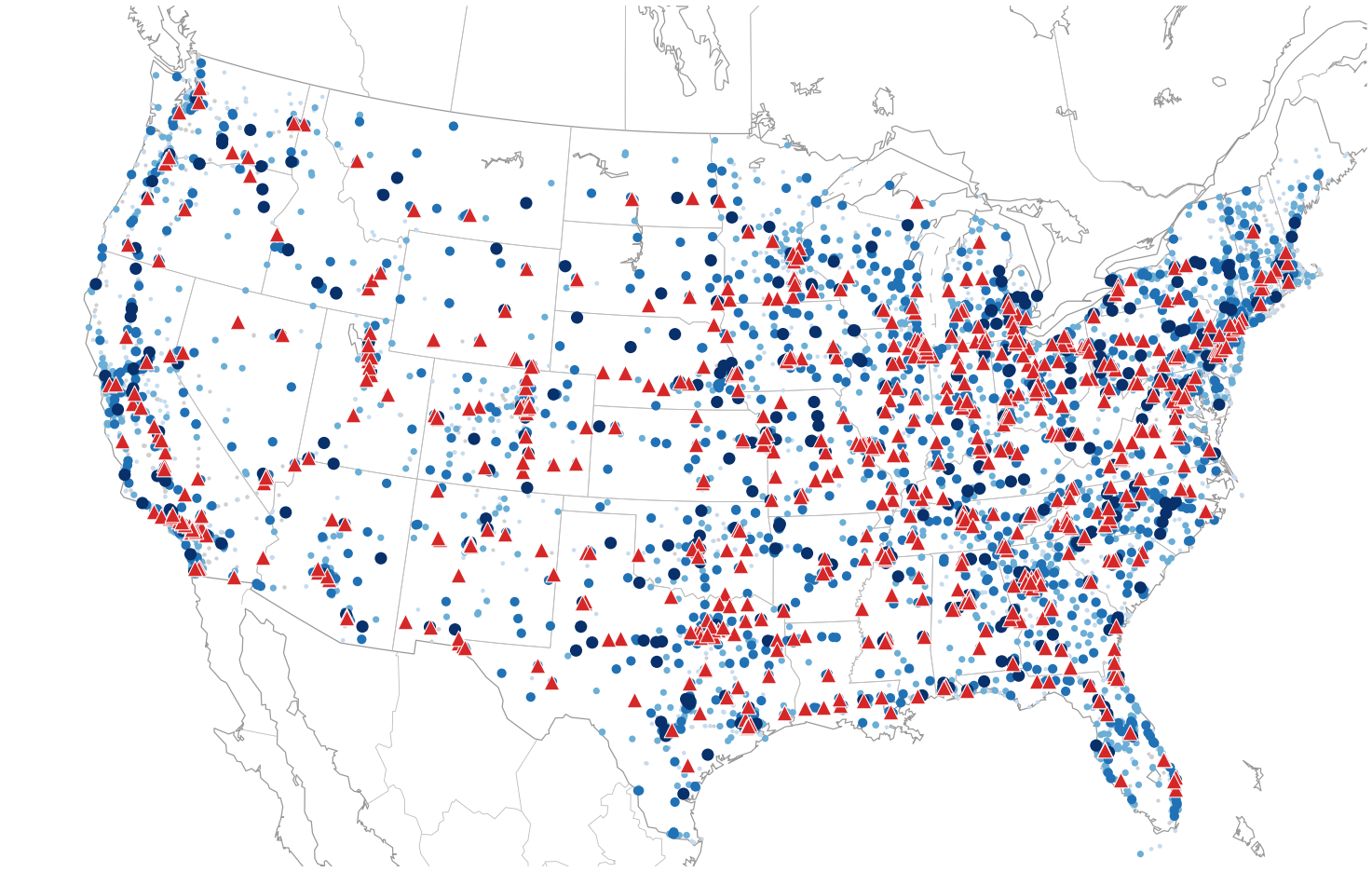}\caption{DCFC 50~kW}\end{subfigure}\hfill
\begin{subfigure}{0.49\textwidth}\includegraphics[width=\linewidth]{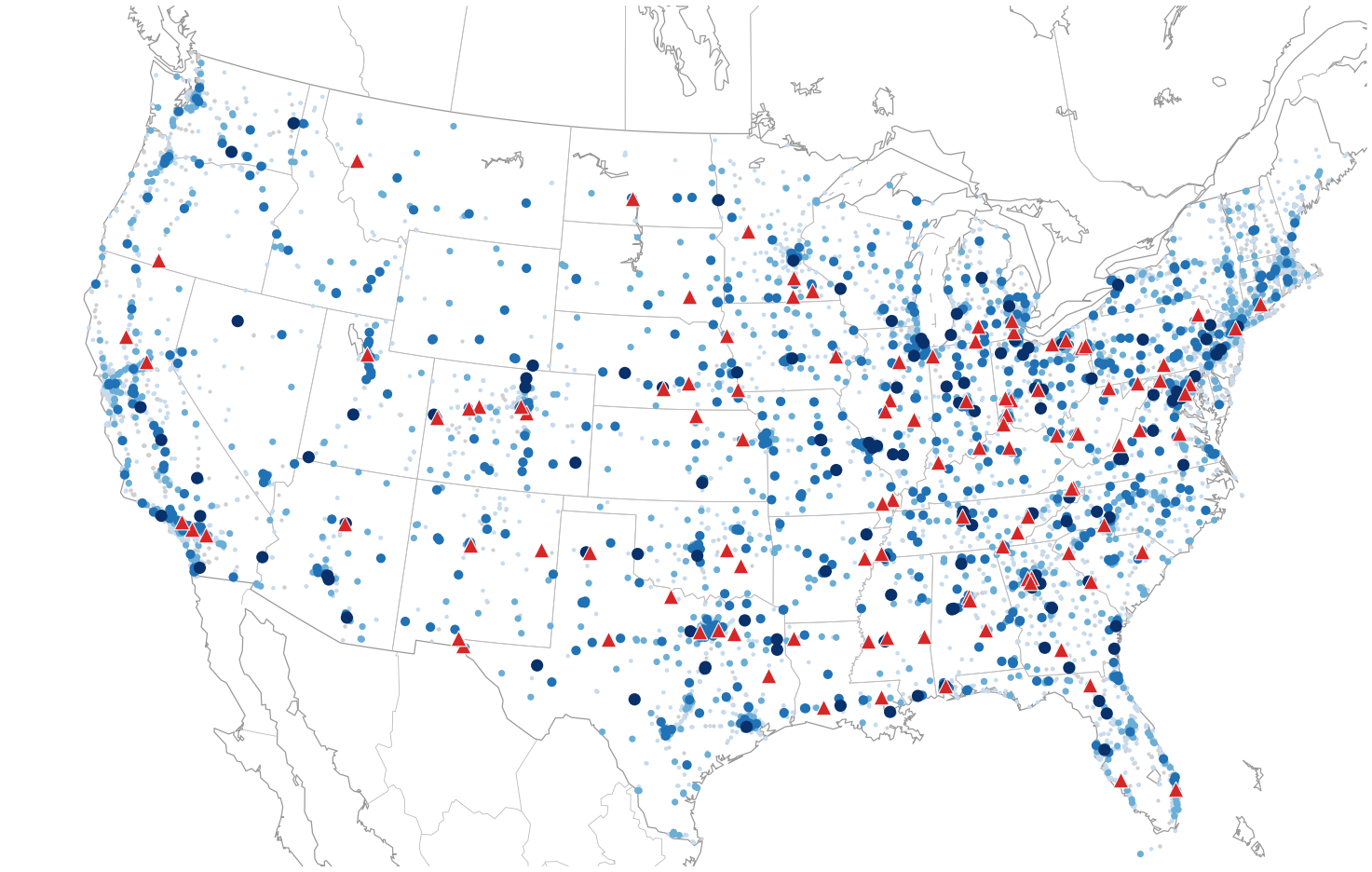}\caption{DCFC 350~kW}\end{subfigure}\\[4pt]
\begin{subfigure}{0.49\textwidth}\includegraphics[width=\linewidth]{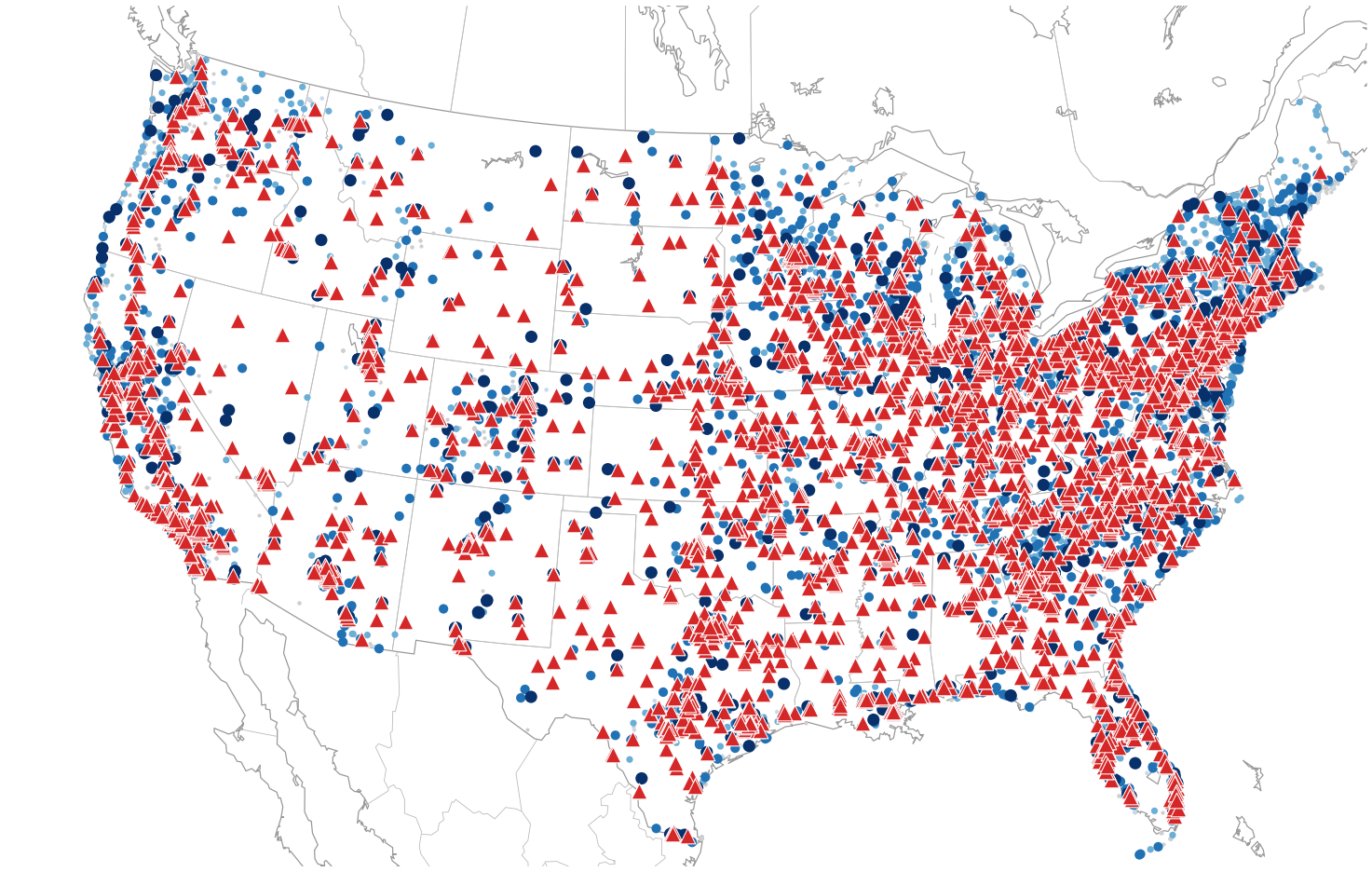}\caption{L2 5~kW}\end{subfigure}\hfill
\begin{subfigure}{0.49\textwidth}\includegraphics[width=\linewidth]{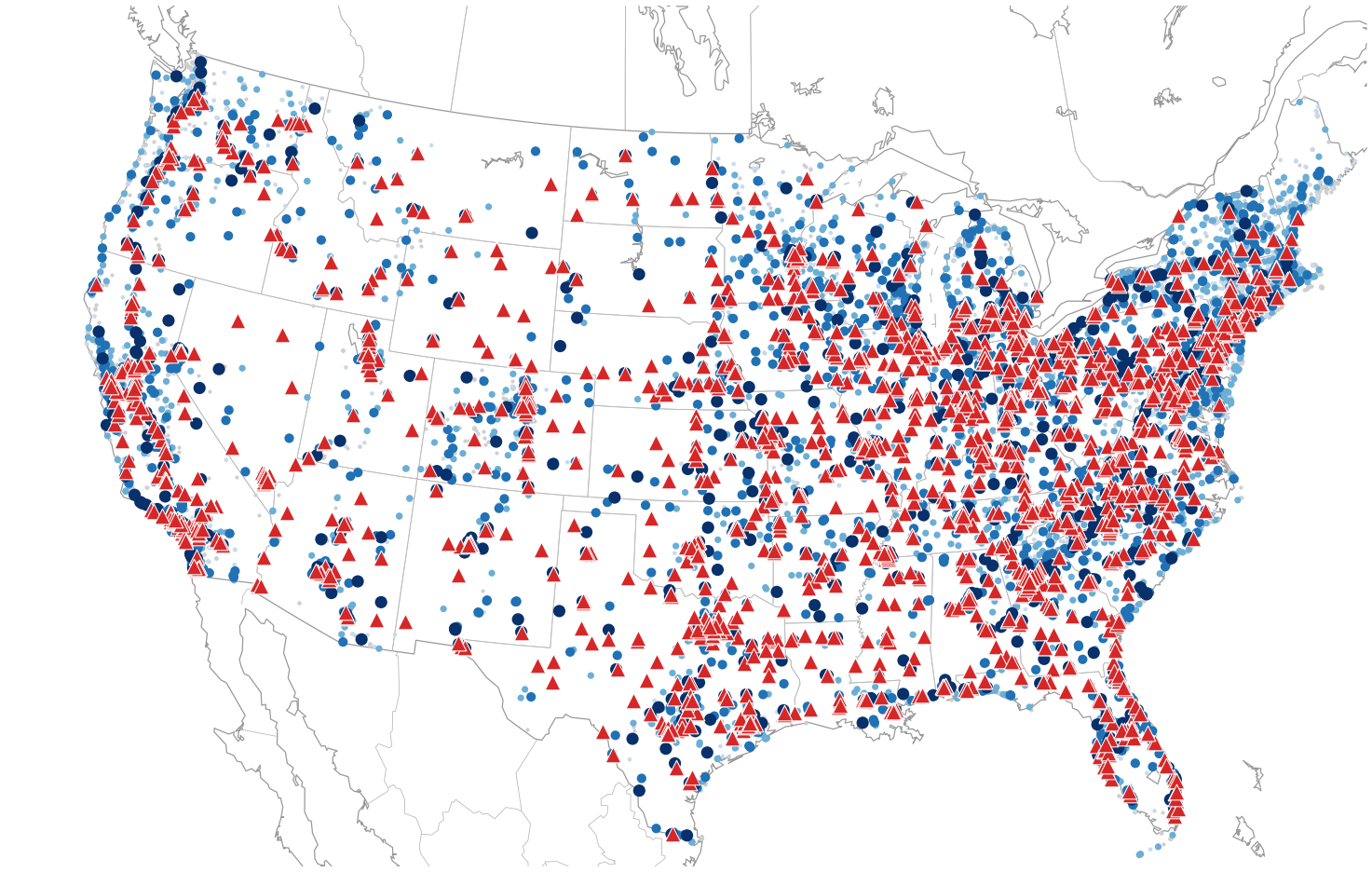}\caption{L2 15~kW}\end{subfigure}\\[2pt]
\includegraphics[width=\textwidth]{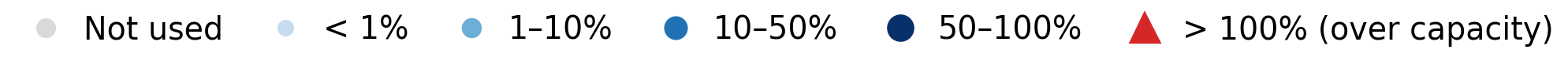}
\caption{Monthly plug utilization of each station in the uncoordinated national runs (DAG + FIFO, three candidate stations per link endpoint; charging hours divided by plug-hours over 30 days). Red triangles: stations whose demand exceeds their capacity (681, 112, 4,440, and 2,224 stations; \cref{tab:cong}). Grey: stations not used. Map lines delineate study areas and do not necessarily depict accepted national boundaries.}\label{fig:util}
\end{figure}

This subsection shows what capacity-blind planning produces: each trip takes the stops that are best for it alone, and plugs are then shared first-come, first-served. This is the starting point of the two re-planning methods of \cref{sec:congestion} and the outcome of large-scale models that do not represent plug capacity; \cref{sec:coordination} reports the coordinated outcome, the other end of the range of traveler behavior we model. At the level of the system, capacity is ample (\cref{tab:cong}): over the month, the plugs of the DCFC network are busy 14\% of the time at 50~kW and 2\% at 350~kW, and the median station used is busy 2\% and 0.2\% of the time. The load is, however, extremely uneven (\cref{fig:util}). At 50~kW, 681 stations, 7\% of those used, receive more charging hours during the month than their plugs can deliver, and they host 63\% of all charging sessions; the 100 stations with the most waiting account for 92\% of it.

\begin{table}[!htb]
\centering
\begin{threeparttable}
\caption{Plug load in the national scenarios without coordination.}\label{tab:cong}
\footnotesize
\begin{tabular}{l r r r r}
\toprule
 & DCFC 50 kW & DCFC 350 kW & L2 5 kW & L2 15 kW \\
\midrule
Charging sessions (million) & 8.07 & 8.16 & 7.71 & 7.72 \\
Stations used & 9,684 & 8,764 & 18,451 & 17,931 \\
Plug utilization, system (\%) & 14.1 & 2.1 & 53.8 & 17.9 \\
Plug utilization, median station (\%) & 2.3 & 0.2 & 16.0 & 5.1 \\
Stations over capacity & 681 & 112 & 4,440 & 2,224 \\
\quad share of stations used (\%) & 7.0 & 1.3 & 24.1 & 12.4 \\
\quad share of sessions (\%) & 62.7 & 24.9 & 92.3 & 83.3 \\
Demand above capacity (\% of charging hours) & 43.5 & 11.3 & 82.3 & 68.6 \\
Trips that wait $>$1 h (\%) & 72.8 & 47.1 & 89.8 & 83.8 \\
\bottomrule
\end{tabular}
\begin{tablenotes}[flushleft]\footnotesize
\item \textit{Note:} DAG + FIFO after the two-hour best-response stage, three candidate stations per link endpoint. Utilization: charging hours divided by plug-hours over the month; a station is over capacity when its utilization exceeds 100\%.
\end{tablenotes}
\end{threeparttable}
\end{table}

In the cases we could test with reservations (\cref{sec:coordination}), the overload comes from the concentration of trips on the same stations rather than from a shortage of plugs in the region. Every trip chooses among the nearest station and up to two others within 15 miles of where it stops and, without information about queues, picks the one that minimizes its own detour and charging time, so trips on the same interstate choose the same few stations; with one candidate per stop (run E of \cref{tab:ablation}) the overload is the same, so the concentration comes mainly from the nearest stations themselves. At 50~kW, 68\% of the over-capacity stations have one or two plugs, against 43\% of all stations used, and 91\% of them have a station with at least four plugs within 15 miles that is busy less than half the time. Adding 4,556 plugs at these stations, 6.7\% of the national DCFC plug count, would cover their monthly demand, though not its peaks within the day and week. Faster charging shortens sessions sevenfold, so at 350~kW only 112 stations are over capacity; on the L2 networks, where a session lasts 2.8--8.5 hours, 12--24\% of the stations used are over capacity and host 83--92\% of the sessions.

Without coordination, queues at over-capacity stations keep growing, and some trips wait hours or even days. These waits measure demand that the chosen stations cannot serve rather than delays a traveler would accept, which is why \cref{tab:econ} excludes waiting and \cref{tab:cong} describes congestion mainly through station counts. Waiting is close to all-or-nothing: at 50 and 350~kW only 2.5\% and 3.5\% of the trips wait between 15 minutes and an hour, so travelers who left a queue after a modest wait would mostly be those who otherwise wait more than an hour. The overload appears at low demand: when we reran the 50~kW scenario on random subsamples of the trips, 9\% of the trips already waited (3\% more than an hour) in a 1\% subsample (about 0.01\% of national demand), and 35\% waited more than an hour in a 10\% sample.

\subsection{Coordination through plug reservations}\label{sec:coordination}
At national scale, the two re-planning heuristics did not change this picture within their limits. With three candidate stations, best response re-planned only 50--250 of the 1.3--2.0 million waiting trips in its two hours, so its national result is effectively the uncoordinated one rather than a converged solution; its quality is better judged on the corridor instances of \cref{sec:benchmark}. The station-penalty heuristic, run for about 25 hours per scenario, did not reduce total waiting (\hyperref[sec:supp]{Supplementary material}, \cref{sec:S2}). We therefore ran sequential reservation (\cref{sec:congestion}) with nearest and with capacity-aware candidate sets on all 2.2 million trips at 350~kW, and on the random 10\% sample of \cref{sec:congestionresults} at 50 and 350~kW, each paired with DAG + FIFO on the same trips. At 50~kW the pass over all trips was too slow: after 14 hours it had planned about 9\% of them, because sessions are seven times longer and far fewer trips can take the fast path. The 10\% sample represents about 0.1\% of national demand, a level at which 35\% of uncoordinated trips already wait more than an hour at 50~kW. \Cref{tab:coord} reports the results.

\begin{table}[!htb]
\centering
\begin{threeparttable}
\caption{Coordination through plug reservations.}\label{tab:coord}
\footnotesize
\setlength{\tabcolsep}{3pt}
\begin{tabular}{l r r r r r r r}
\toprule
Method & Wait & Detour & Stops & Stations & Tier 1, with & Cheaper & Run time \\
 & $>$1 h (\%) & (min/trip) & per trip & used & wait (\%) & EV (\%) & (h) \\
\midrule
\multicolumn{8}{l}{\emph{10\% of trips (220,809), 50~kW}}\\
DAG + FIFO & 35.4 & 8.00 & 3.65 & 7,936 & --- & 4.1 & 2.1 \\
Res., nearest & 0.0 & 8.47 & 3.67 & 8,456 & $-$145.3 & 4.5 & 8.9 \\
Res., capacity-aware & 0.0 & 8.58 & 3.68 & 8,196 & $-$145.5 & 4.5 & 7.7 \\
\midrule
\multicolumn{8}{l}{\emph{10\% of trips (220,809), 350~kW}}\\
DAG + FIFO & 1.9 & 6.88 & 3.69 & 6,864 & +17.2 & 81.2 & 2.1 \\
Res., nearest & 0.0 & 6.94 & 3.70 & 7,097 & +21.4 & 84.3 & 4.3 \\
Res., capacity-aware & 0.0 & 6.98 & 3.70 & 6,961 & +21.4 & 84.3 & 4.3 \\
\midrule
\multicolumn{8}{l}{\emph{All feasible trips (2,210,571), 350~kW}}\\
DAG + FIFO & 47.1 & 6.63 & 3.69 & 8,764 & --- & 43.3 & 20 \\
Res., nearest & 0.0 & 7.33 & 3.72 & 9,474 & +20.8 & 83.6 & 67.1 \\
Res., capacity-aware & 0.0 & 7.44 & 3.73 & 9,079 & +20.7 & 83.5 & 61.6 \\
\bottomrule
\end{tabular}
\begin{tablenotes}[flushleft]\footnotesize
\item \textit{Note:} DAG + FIFO: uncoordinated plans queued FIFO (on all trips, after the two-hour best-response stage, which changed less than 0.5\% of the cost). Res.: sequential reservation with the nearest or the capacity-aware candidate set. All plans are scored with the same FIFO simulation with delay propagation. Wait $>$1 h: share of trips that wait more than one hour; under reservation, 4 trips (10\%, 50~kW) and 7 trips (all, 350~kW) wait more than an hour, which rounds to 0.0. The 10\% sample is drawn from all 2.77 million trips before the filters of \cref{sec:design}. Tier~1: reduction of the total EV operating cost \cref{eq:tier1}, including waiting, relative to the ICE cost at \$0.15 per kWh; negative values mean the EV is more expensive; ---: not reported where long waits under unlimited patience dominate the cost (more than a third of the trips wait over an hour). Cheaper EV: share of trips whose EV Tier~1 cost is below the ICE cost. Run time: wall clock on a shared 128-core node, end to end.
\end{tablenotes}
\end{threeparttable}
\end{table}

Four results stand out. First, planning against reservations removes the waiting almost entirely, on the same network and for the same trips. On all trips at 350~kW, uncoordinated scheduling leaves 47\% of the trips waiting more than an hour; with reservations seven of the 2.2 million trips wait more than an hour, the mean wait is 0.07 minutes, and 95\% of the trips wait less than 0.4 minutes. On the 10\% sample at 50~kW, where 35\% of uncoordinated trips wait more than an hour, the mean wait falls to about a second. The remaining waits are short: on all trips at 350~kW, 99\% of the trips wait less than two minutes, and the longest wait is under three hours. Some are short waits that the plans accept (\cref{sec:congestion}); most are likely because the reservation pass and the evaluation compute arrival times separately; the FIFO evaluation can also let an earlier arrival take a plug booked for a later trip, and we did not trace the few long waits (\hyperref[sec:supp]{Supplementary material}, \cref{sec:S3}).

Second, coordination costs little. Trips avoid booked plugs by stopping at another nearby station or at another point of the route: at 350~kW on all trips, the mean detour grows by 0.7 minutes per trip, stops by 0.03 per trip, and the number of stations used by 8\%. Because waiting nearly disappears, the cost comparison with waiting included becomes the comparison of \cref{tab:econ}, which excludes it: at 350~kW the EV trip is 20.8\% cheaper than the ICE trip in total (22.0\% in \cref{tab:econ}, whose plans have shorter detours) and cheaper for 84\% of the trips, and at 50~kW it is 145\% more expensive, as in the sensitivity baseline. Per-trip detours and costs in the 10\% sample and the 1\% subsample differ slightly from those of the full data set in \cref{tab:case,tab:econ} (by up to 0.8 minutes of detour and two points of Tier~1; the 1\% subsample also includes an effective best-response stage, which adds detour). Comparisons are therefore made within one sample.

Third, capacity-aware candidate sets add nothing once trips reserve. They are no better than the nearest candidates, are slightly worse on most measures though never slower, and use fewer stations, because reserving two of the three slots for large stations removes nearby small stations that are often free when a trip needs them. What trips need is a free plug at the time they arrive, and the reservation calendar finds one among the nearest stations. We did not test capacity-aware candidates without reservation, where they may matter more.

Fourth, the pass is slow but feasible at national scale. On all trips at 350~kW the reservation pass itself took 20--23 hours, and the run 62--67 hours end to end, about 10 hours of which wrote the link-level SoC trajectories; the rest went to loading the data, solving the capacity-free DP, and evaluating the schedules on a shared node. About 44\% of the trips took the fast path. At 50~kW, overloaded stations hold long sequences of reserved sessions, more trips need the time-dependent DP, and the pass over all trips did not finish in reasonable time; on the 10\% sample it took about four times as long as DAG + FIFO. The comparison with best response is not like for like: best response starts from a FIFO state in which most trips already wait and had a two-hour budget, whereas the reservation pass starts from the capacity-free plans and ran to completion. A best-response stage started from the reservation plans, or given a comparable budget, may do better than the one we ran, although in the one run in which we logged such a round it found no trip that could lower its own cost (\hyperref[sec:supp]{Supplementary material}, \cref{sec:S3}).

These results change how the congestion results of \cref{sec:congestionresults} should be read. In the cases we could run, the overload comes from uncoordinated use of the same stations rather than from a shortage of plugs: when each trip knows which plug-times are taken, the existing DCFC network serves all 2.2 million trips at 350~kW, and the 10\% sample (about 0.1\% of national demand) at 50~kW, with almost no waiting. We have not shown this for all trips at 50~kW or for L2. Reservation is also an idealization: a booking service that every traveler uses and follows, with on-time departures, no local charging demand at the same plugs, and planning in order of departure, whose outcome depends on that order. The results are therefore an optimistic case for booking, which assumes full compliance; \cref{sec:discussion} discusses what this implies.

\subsection{Sensitivity}\label{sec:sens}
\begin{table}[!tb]
\centering
\begin{threeparttable}
\caption{Sensitivity to EV range and charger power on the 1\% subsample.}\label{tab:sens}
\footnotesize
\begin{tabular}{l r r r r r r r}
\toprule
Level & Tier 1 (\%) & Tier 2 (\%) & Stops & Detour & Charge & $\Delta d$ (\%) & Wait (\%) \\
\midrule
\multicolumn{8}{l}{\emph{EV range $\underline{R}$--$\overline{R}$ (miles)}}\\
50--100 & $-$247.5 (1.3) & $-$83.8 (0.4) & 8.83 & 24.5 & 18.2 & 4.66 & 0.57 \\ 
100--200 & $-$185.5 (0.5) & $-$62.0 (0.2) & 5.69 & 13.8 & 33.4 & 1.99 & 0.29 \\ 
150--350$^*$ & $-$145.7 (0.6) & $-$48.3 (0.2) & 3.68 & 8.5 & 51.6 & 1.07 & 0.13 \\ 
200--300 & $-$144.5 (0.7) & $-$47.9 (0.2) & 3.51 & 8.0 & 54.2 & 1.01 & 0.13 \\ 
300--400 & $-$115.6 (0.4) & $-$38.1 (0.1) & 2.53 & 5.7 & 72.6 & 0.66 & 0.06 \\ 
400--500 & $-$93.1 (0.5) & $-$30.5 (0.2) & 2.00 & 4.8 & 88.2 & 0.50 & 0.04 \\ 
\midrule
\multicolumn{8}{l}{\emph{Charger power $W$ (kW)}}\\
50$^*$ & $-$145.4 (0.5) & $-$48.2 (0.2) & 3.68 & 8.5 & 51.5 & 1.07 & 0.12 \\ 
100 & $-$48.8 (0.5) & $-$15.6 (0.2) & 3.69 & 7.9 & 25.8 & 1.00 & 0.15 \\ 
250 & +9.7 (0.2) & +4.1 (0.1) & 3.71 & 7.4 & 10.4 & 0.93 & 0.15 \\ 
350 & +21.0 (0.2) & +8.0 (0.1) & 3.73 & 7.3 & 7.4 & 0.92 & 0.13 \\ 
500 & +29.4 (0.2) & +10.8 (0.1) & 3.73 & 7.2 & 5.2 & 0.90 & 0.13 \\ 
\bottomrule
\end{tabular}
\begin{tablenotes}[flushleft]\footnotesize
\item \textit{Note:} DCFC 50~kW unless varied; $^*$ baseline. Tier~1 and Tier~2: percentage reduction of the EV cost relative to the ICE cost, including waiting, as the mean over ten seeds with the standard deviation across seeds in parentheses. Stops: per trip; detour: minutes per trip; charge: minutes per stop; $\Delta d$: added distance (\%); wait: share of trips with a nonzero wait (\%). Values are over the trips that need to charge, a set that changes with range; waiting is handled by best response, which is effective at this size. All sweeps are reported in \cref{tab:S1,tab:S2} of the \hyperref[sec:supp]{Supplementary material}.
\end{tablenotes}
\end{threeparttable}
\end{table}

\Cref{tab:sens} reports the sweeps of EV range and charger power on the 1\% subsample (about 22,000 trips that need to charge, ten seeds per level); \cref{tab:S1,tab:S2} in the \hyperref[sec:supp]{Supplementary material} report all sweeps. With best-response coordination, which is effective at this size, congestion is negligible: fewer than 1\% of trips wait in any configuration, even after removing 90\% of the plugs (\hyperref[sec:supp]{Supplementary material}, \cref{tab:S2}), so the sweeps isolate the effect of each parameter on routing and cost. Without coordination, 9\% of the trips of a random sample of the same size have a nonzero wait (\cref{sec:congestionresults}). The standard deviations across seeds, shown in parentheses in the tables, are at most 1.3 percentage points. Each sweep was run with its own baseline replicates; their Tier~1 values range from $-145.4$\% to $-145.8$\%, which indicates the run-to-run variation, and differences larger than about one point within a sweep are not noise.

Charger power has the largest effect: the Tier~1 comparison moves from $-145$\% at 50~kW to $+10$\% at 250~kW and $+29$\% at 500~kW. EV range has the next largest effect, through the number of stops: at 50--100 miles of range trips need 8.8 stops and at 400--500 miles 2.0. Energy use, SoC limits, and the starting SoC change the cost by 10--90 points through the energy charged (\hyperref[sec:supp]{Supplementary material}, \cref{tab:S1}). The station-choice and network parameters (\hyperref[sec:supp]{Supplementary material}, \cref{tab:S2}), namely the candidate-set size $K$, detour circuity, and detour speed, change Tier~1 by at most four points from the baseline. Costs and routing are nearly the same for $K = 1$, 3, and 5 (only the small share of trips that wait changes), which supports the three candidates used in the national runs; beyond the nearest station, the next candidates are usually farther than the detour is worth when there is no queue. Raising the circuity factor to 1.4 changes Tier~1 by about three points, so the $L_1$ detours are not a material simplification at 50~kW.

\subsection{Computation}\label{sec:runtime}
The national runs were made on 128-core nodes of the Laboratory Computing Resource Center at Argonne National Laboratory, with several scenarios sharing a node. Loading the link-level trajectories and solving the capacity-free DP for 2.2 million trips took 6--9 hours of wall-clock time per scenario in this shared setting, with no commercial solver. The best-response stage used its full two-hour limit, and writing the link-level SoC trajectories took most of the remaining time, for about 20 hours end to end. The penalty heuristic with three candidate stations took about 25 hours per scenario. Sequential reservation over all trips at 350~kW took 62--67 hours end to end (\cref{sec:coordination}). These times are dominated by data handling and by the congestion stage; the per-vehicle DP itself is linear in the number of DAG arcs and parallel across vehicles.
\section{Discussion and conclusions}\label{sec:discussion}
This paper scheduled the charging stops of 2.2 million electrified long-distance trips through the U.S.\ public charging network with a per-vehicle dynamic program that is exact under the charge-to-target policy and needs no commercial solver. Under the model's assumptions (every DCFC plug usable by every EV, the nearest station always a candidate), the existing DCFC network reaches 98.5\% of the trips that need to charge (at least 97.8\% if the 15,711 trips with an overly long link count as unserved), with fewer than four stops and about seven minutes of detour per trip; gaps in coverage of the kind identified by \citet{hanig_2025} matter for particular corridors but do not make many trips infeasible. What separates EVs from ICE vehicles on long trips is the time spent charging and the price of the energy bought en route. At \$0.15 per kWh, the EV trip is cheaper only for travelers who value their charging time below about \$11 per hour at 50~kW, or \$111 per hour at 350~kW, and at 250~kW the advantage at baseline values shrinks to about 10\% (\cref{sec:sens}); at public fast-charging prices of \$0.35--0.45 per kWh, it is more expensive at either power for practically any value of time.

How the problem is modeled matters as much as the parameters. Three choices that are easy to make in a large model, namely rejecting trips that start above the charging target, discretizing SoC, and omitting the value of charging time, each changed feasibility, stops, or the sign of the cost comparison by a large factor. The SoC-bucket error is general: any DP that rounds each link's energy up to a whole bucket overstates energy on networks with many short links, whereas the opportunity DAG needs no buckets because the SoC after every stop is known in advance.

Where we could test reservations (all trips at 350~kW and a 10\% sample at 50~kW), congestion is mainly a coordination problem. Without coordination, several hundred mostly small stations are overloaded while larger stations nearby are busy less than half the time, and waiting appears at a hundredth of a percent of national demand. Sequential reservation removes almost all of this waiting for less than a minute of extra detour per trip, on all trips at 350~kW and on a 10\% sample at 50~kW: the same network serves the same trips once each trip knows when a plug is free. The uncoordinated and reservation runs are the two ends of a range: travelers with no information about plugs, and travelers who all book and keep their bookings. The re-planning heuristics did not achieve this at national scale. Best response, reliable on small corridors, re-planned fewer than 600 trips within its budget, and the station-penalty heuristic re-prices all congested stations at once, so the re-planned trips move together to stations that may or may not have room; a penalty scheme for capacitated charging needs prices that accumulate over iterations, as dual prices do, or an acceptance test against the evaluated cost. Reservation is 6\% from the best known solution on small corridors, so its national schedules are not optimal, but it is the only method we tested that resolved the congestion of a national instance, at 350~kW. The aggregate figures also bound the future: about 1\% of national long-distance demand uses 14\% of the DCFC plug-hours at 50~kW, so, if all long-distance trips resembled this sample, serving them would need at least 14 times today's DCFC plug capacity at 50~kW, or at least twice at 350~kW, even at full utilization, before any local demand. At the 1\% level and 350~kW, a booking service is a lever that needs no new plugs; at 50~kW without one, about 4,600 additional plugs would be needed to match monthly demand, and even they would not remove waiting at peak times.

These results rest on simplifications. Every stop charges to 80\% at constant power; a higher target mainly reduces the number of stops, charging every final stop to $\overline{S}$ lets trips arrive with more charge than they need, which overstates the EV's charging time slightly (\cref{sec:econresults}), and the opportunity DAG extends exactly to any finite set of targets, whereas a continuous charge amount or a nonlinear charging curve \citep{montoya_2017} would require labels with SoC again. Detours use $L_1$ distances, to which the results at 50~kW are not sensitive. Uncoordinated travelers do not react to queues, whereas in the reservation runs every traveler books, keeps the booking, and departs on time, with no local demand at the same plugs; real behavior lies between the two, and we did not test an intermediate policy in which travelers react to observed queues, and the reservation results are an optimistic case for booking. The time-dependent DP keeps one label per node, so it is exact without waiting and heuristic with it. All plugs of a scenario share one power, every DCFC plug is assumed usable by every EV and in service, electricity is priced at a single rate unless stated, energy use and range are drawn independently, one cost per plug is used for all plug types, meal and rest breaks enter only through the bound of \cref{sec:econresults}, and home, workplace, and destination charging enter only through the starting SoC, and overnight stops on trips of 500 miles or more, where travelers could charge while they sleep, are not modeled. Finally, the full-sample runs represent one demand level, about 1\% of national long-distance demand, taken from POLARIS without separate validation of long-distance trip frequencies against survey data, and coordinated runs at national scale were made only at 0.1\% (50~kW) and at 0.1\% and 1\% (350~kW) of national long-distance demand, besides the 0.01\% subsample sweeps with best response. Four extensions follow directly: faster reservation passes for low-power, heavily loaded networks (for example, by region and time window); tests of planning order, partial adoption, no-shows, and background local demand; behavioral models of how travelers use reservations and react to queues; and coupling the scheduler with plug-expansion decisions for the demand levels at which reservations alone no longer suffice.

\section*{CRediT authorship contribution statement}
\textbf{Taner Cokyasar:} Conceptualization, Methodology, Software, Formal analysis, Investigation, Data curation, Visualization, Writing -- original draft, Writing -- review \& editing.

\section*{Declaration of competing interest}
The author declares that he has no known competing financial interests or personal relationships that could have appeared to influence the work reported in this paper.

\section*{Funding}
This work was supported by the U.S. Department of Energy, Office of Science [contract number DE-AC02-06CH11357], and sponsored by the DOE Transportation Technologies Office under the Pathways to Affordable, Convenient, and Efficient Regional Mobility initiative of the Energy Efficient Mobility Systems Program. A DOE program manager helped establish the project concept, advance its implementation, and provide guidance (see Acknowledgments).

\section*{Data and code availability}
The results needed to interpret the findings are reported in the tables of this article. The scheduling code is part of the POLARIS tools, and the trip trajectories come from the POLARIS model; both are licensed by Argonne National Laboratory and cannot be redistributed by the author.

\section*{Acknowledgments}
This material is based on work supported by the U.S. Department of Energy (DOE), Office of Science, under contract number DE-AC02-06CH11357. This report and the work described were sponsored by the DOE Transportation Technologies Office (TTO) under the Pathways to Affordable, Convenient, and Efficient Regional Mobility, an initiative of the Energy Efficient Mobility Systems (EEMS) Program. Erin Boyd, a DOE Office of Critical Minerals and Energy Innovation (CMEI) manager, played an important role in establishing the project concept, advancing implementation, and providing guidance.

\section*{Declaration of generative AI and AI-assisted technologies in the manuscript preparation process}
During the preparation of this work the author used Claude (Anthropic) to assist with code development, data-analysis scripts, drafting, and language editing, and ChatGPT (OpenAI), Gemini (Google), Grok (xAI), Perplexity, and Mistral to obtain critical feedback on drafts. After using these tools, the author reviewed and edited the content as needed and takes full responsibility for the content of the published article.

\bibliographystyle{elsarticle-harv}

\clearpage
\phantomsection\label{sec:supp}
\section*{Supplementary material}
\setcounter{section}{0}\setcounter{table}{0}
\renewcommand{\thesection}{S\arabic{section}}\renewcommand{\theHsection}{S\arabic{section}}
\renewcommand{\thetable}{S\arabic{table}}\renewcommand{\theHtable}{S\arabic{table}}

\section{Corridor instances}\label{sec:S1}
We generated 25 corridor instances small enough for the exact model~\crefrange{eq:obj}{eq:dom}. Each has a 600-mile highway split into links of 2--8 miles, 18 single-plug stations placed 0.2--2.5 miles off the highway, and $n \in \{5, 10, 15, 20, 30\}$ trips (five random instances each) of 220--400 miles that depart within the same hour. Each trip has a range of 180--300 miles, energy use $\sim \mathrm{Triangular}(250, 250, 450)$~Wh/mile, a starting SoC $\sim \mathrm{Triangular}(60, 100, 100)\%$, and, as candidates of each link endpoint, those of its $K = 2$ nearest stations that lie within 4 miles (unlike (A3), a nearest station farther away is not a candidate); chargers deliver 50~kW and $C_\tau$, $C_\chi$, and $C_\omega$ are \$70, \$50, and \$50 per hour. These instances are deliberately congested: many trips need the same few single-plug stations within a short window, a harder setting for coordination than the national data, whose stations have more plugs and whose departures spread over a month. Energy use is drawn from a narrower range than in the national runs; the instances are meant to test coordination, not range.

\section{Station-penalty stage in the stepwise comparison}\label{sec:S2}
The last change of the stepwise comparison, run E, which replaces station penalties with FIFO best response, shows how much the congestion stage matters at national scale, and in which direction. In all of B--D the penalty stage stopped after a single iteration, because that iteration alone outlasted the stage's one-hour time limit in these runs (the case-study runs of \cref{sec:coordination} used a longer limit and took about 25 hours); it re-planned about two million trips at once against penalties on the 2,700--3,700 stations whose occupancy exceeded their plugs in at least one 15-minute bin, a much weaker condition than the monthly over-capacity count of \cref{tab:ablation}. In run E, best response re-planned only 560 trips in its two hours, so run E is in effect the uncoordinated DAG + FIFO outcome. Compared with it, the single penalty iteration of run D spread the load (398 instead of 682 stations over capacity, with $K = 1$), and lowered the share of trips waiting more than an hour from 73\% to 61\%, at the cost of 13 more minutes of detour per trip. We do not read this as evidence that the heuristic coordinates well: one simultaneous re-pricing happened to spread the load in this configuration, and, as \cref{sec:coordination} notes, the same heuristic with three candidate stations per link endpoint did not reduce waiting. In the national case runs with three candidates, the penalty stage raised the share of trips waiting more than an hour from 72.8\% to 75.1\% at 50~kW and lowered it from 47.1\% to 43.1\% at 350~kW, while adding about 5 and 3.5 minutes of detour per trip; total waiting rose in both.

\section{Waits under sequential reservation}\label{sec:S3}
With sequential reservation, the mean evaluated wait is 0.07 minutes per trip on all trips at 350~kW and about a second on the 10\% sample at 50~kW (\cref{sec:coordination}). About half of the trips record a nonzero wait, but most of these waits last seconds (95\% of all trips wait less than 0.4 minutes); their likely source is that the reservation pass and the evaluation simulation compute arrival and session times separately, so a booked session can start a few seconds late in the evaluation. A few trips wait longer: on all trips at 350~kW, 0.02\% wait more than ten minutes and seven trips more than an hour, the longest 2.8 hours, and these seven account for 0.4\% of the total waiting; on the 10\% sample at 50~kW, four trips wait more than an hour. We did not reconcile the two clocks, so we report the evaluated waits as they are. The few long waits may also arise because the FIFO evaluation lets an earlier arrival take a plug booked for a later trip; we did not trace them. In the one run in which we logged it (50~kW, nearest candidates, 10\% of trips), a best-response round after the pass found no trip that could lower its own cost.

\section{Sensitivity sweeps}\label{sec:S4}
\Cref{tab:S1,tab:S2} report all sweeps on the 1\% subsample (\cref{sec:sens}).

\begin{table}[!htb]
\centering
\begin{threeparttable}
\caption{Sensitivity to vehicle and charging parameters on the 1\% subsample.}\label{tab:S1}
\footnotesize
\begin{tabular}{l r r r r r r r}
\toprule
Level & Tier 1 (\%) & Tier 2 (\%) & Stops & Detour & Charge & $\Delta d$ (\%) & Wait (\%) \\
\midrule
\multicolumn{8}{l}{\emph{EV range $\underline{R}$--$\overline{R}$ (miles)}}\\
50--100 & $-$247.5 (1.3) & $-$83.8 (0.4) & 8.83 & 24.5 & 18.2 & 4.66 & 0.57 \\ 
100--200 & $-$185.5 (0.5) & $-$62.0 (0.2) & 5.69 & 13.8 & 33.4 & 1.99 & 0.29 \\ 
150--350$^*$ & $-$145.7 (0.6) & $-$48.3 (0.2) & 3.68 & 8.5 & 51.6 & 1.07 & 0.13 \\ 
200--300 & $-$144.5 (0.7) & $-$47.9 (0.2) & 3.51 & 8.0 & 54.2 & 1.01 & 0.13 \\ 
300--400 & $-$115.6 (0.4) & $-$38.1 (0.1) & 2.53 & 5.7 & 72.6 & 0.66 & 0.06 \\ 
400--500 & $-$93.1 (0.5) & $-$30.5 (0.2) & 2.00 & 4.8 & 88.2 & 0.50 & 0.04 \\ 
\midrule
\multicolumn{8}{l}{\emph{Energy-use upper bound $\overline{\rho}$ (Wh/mile)}}\\
300 & $-$61.4 (0.1) & $-$19.9 (0.0) & 3.69 & 8.2 & 33.0 & 1.03 & 0.12 \\ 
400 & $-$80.0 (0.2) & $-$26.2 (0.1) & 3.68 & 8.2 & 37.2 & 1.04 & 0.14 \\ 
500 & $-$98.9 (0.4) & $-$32.5 (0.1) & 3.69 & 8.3 & 41.2 & 1.05 & 0.16 \\ 
600 & $-$117.5 (0.4) & $-$38.8 (0.1) & 3.68 & 8.4 & 45.4 & 1.07 & 0.14 \\ 
750$^*$ & $-$145.8 (0.6) & $-$48.3 (0.2) & 3.68 & 8.5 & 51.6 & 1.07 & 0.12 \\ 
\midrule
\multicolumn{8}{l}{\emph{Charger power $W$ (kW)}}\\
50$^*$ & $-$145.4 (0.5) & $-$48.2 (0.2) & 3.68 & 8.5 & 51.5 & 1.07 & 0.12 \\ 
100 & $-$48.8 (0.5) & $-$15.6 (0.2) & 3.69 & 7.9 & 25.8 & 1.00 & 0.15 \\ 
250 & +9.7 (0.2) & +4.1 (0.1) & 3.71 & 7.4 & 10.4 & 0.93 & 0.15 \\ 
350 & +21.0 (0.2) & +8.0 (0.1) & 3.73 & 7.3 & 7.4 & 0.92 & 0.13 \\ 
500 & +29.4 (0.2) & +10.8 (0.1) & 3.73 & 7.2 & 5.2 & 0.90 & 0.13 \\ 
\midrule
\multicolumn{8}{l}{\emph{SoC limits $\underline{S}$--$\overline{S}$ (\%)}}\\
10--100 & $-$132.7 (0.7) & $-$43.8 (0.2) & 2.59 & 5.9 & 73.1 & 0.72 & 0.07 \\ 
20--100 & $-$143.7 (0.6) & $-$47.5 (0.2) & 2.90 & 6.6 & 66.5 & 0.83 & 0.08 \\ 
30--100 & $-$155.7 (0.5) & $-$51.6 (0.2) & 3.28 & 7.4 & 59.8 & 0.96 & 0.10 \\ 
20--80$^*$ & $-$145.7 (0.6) & $-$48.3 (0.2) & 3.68 & 8.5 & 51.5 & 1.08 & 0.13 \\ 
30--80 & $-$159.2 (0.6) & $-$52.9 (0.2) & 4.35 & 10.2 & 44.0 & 1.35 & 0.18 \\ 
\midrule
\multicolumn{8}{l}{\emph{Starting-SoC lower bound $\underline{B}$ (\%)}}\\
60$^*$ & $-$145.8 (0.7) & $-$48.4 (0.2) & 3.68 & 8.5 & 51.6 & 1.07 & 0.14 \\ 
70 & $-$141.7 (0.6) & $-$47.0 (0.2) & 3.67 & 8.4 & 51.4 & 1.06 & 0.12 \\ 
80 & $-$138.1 (0.6) & $-$45.7 (0.2) & 3.64 & 8.4 & 51.3 & 1.05 & 0.14 \\ 
90 & $-$134.8 (1.0) & $-$44.7 (0.3) & 3.63 & 8.4 & 51.2 & 1.03 & 0.13 \\ 
99 & $-$131.9 (0.7) & $-$43.7 (0.2) & 3.61 & 8.4 & 51.2 & 1.02 & 0.12 \\ 
\bottomrule
\end{tabular}
\begin{tablenotes}[flushleft]\footnotesize
\item \textit{Note:} DCFC 50~kW unless varied; $^*$ baseline. Tier~1 and Tier~2: percentage reduction of the EV cost relative to the ICE cost, including waiting, as the mean over ten seeds with the standard deviation across seeds in parentheses. Stops: per trip; detour: minutes per trip; charge: minutes per stop; $\Delta d$: added distance (\%); wait: share of trips with a nonzero wait (\%).
\end{tablenotes}
\end{threeparttable}
\end{table}

\begin{table}[!htb]
\centering
\begin{threeparttable}
\caption{Sensitivity to station-choice and network parameters on the 1\% subsample.}\label{tab:S2}
\footnotesize
\begin{tabular}{l r r r r r r r}
\toprule
Level & Tier 1 (\%) & Tier 2 (\%) & Stops & Detour & Charge & $\Delta d$ (\%) & Wait (\%) \\
\midrule
\multicolumn{8}{l}{\emph{Candidate stations $K$}}\\
1 & $-$145.8 (0.7) & $-$48.3 (0.3) & 3.68 & 8.5 & 51.5 & 1.08 & 0.24 \\ 
3$^*$ & $-$145.8 (0.5) & $-$48.3 (0.2) & 3.68 & 8.5 & 51.6 & 1.07 & 0.13 \\ 
5 & $-$145.5 (0.3) & $-$48.2 (0.1) & 3.68 & 8.5 & 51.5 & 1.08 & 0.13 \\ 
\midrule
\multicolumn{8}{l}{\emph{Detour circuity $\vartheta$}}\\
1$^*$ & $-$145.5 (0.5) & $-$48.2 (0.2) & 3.68 & 8.5 & 51.5 & 1.07 & 0.14 \\ 
1.2 & $-$147.0 (0.5) & $-$48.8 (0.2) & 3.63 & 9.6 & 51.7 & 1.23 & 0.14 \\ 
1.4 & $-$148.9 (0.6) & $-$49.5 (0.2) & 3.58 & 10.7 & 52.2 & 1.37 & 0.15 \\ 
\midrule
\multicolumn{8}{l}{\emph{Detour speed $\bar{s}$ (mph)}}\\
30 & $-$149.5 (0.4) & $-$49.6 (0.1) & 3.68 & 10.9 & 51.6 & 1.04 & 0.19 \\ 
40$^*$ & $-$145.5 (0.7) & $-$48.3 (0.3) & 3.69 & 8.5 & 51.5 & 1.07 & 0.13 \\ 
55 & $-$142.3 (0.6) & $-$47.2 (0.2) & 3.68 & 6.4 & 51.5 & 1.11 & 0.11 \\ 
\midrule
\multicolumn{8}{l}{\emph{Plug supply multiplier}}\\
1$^*$ & $-$145.6 (0.6) & $-$48.3 (0.2) & 3.69 & 8.5 & 51.5 & 1.07 & 0.13 \\ 
0.5 & $-$145.7 (0.4) & $-$48.3 (0.1) & 3.69 & 8.5 & 51.5 & 1.08 & 0.27 \\ 
0.25 & $-$146.3 (0.7) & $-$48.5 (0.2) & 3.69 & 8.7 & 51.5 & 1.10 & 0.47 \\ 
0.1 & $-$146.2 (0.7) & $-$48.5 (0.2) & 3.70 & 8.8 & 51.4 & 1.12 & 0.70 \\ 
\bottomrule
\end{tabular}
\begin{tablenotes}[flushleft]\footnotesize
\item \textit{Note:} DCFC 50~kW; $^*$ baseline. Columns as in \cref{tab:S1}: Tier~1 and Tier~2 are means over ten seeds, with the standard deviation across seeds in parentheses.
\end{tablenotes}
\end{threeparttable}
\end{table}

\end{document}